\documentclass[a4paper,leqno]{article}

\usepackage{hyperref}
\usepackage[T1]{fontenc}
\usepackage[english]{babel}
\usepackage{amsmath}
\input xy \xyoption{all}

\newcounter{szn}
\newcounter{env}[szn]

\renewcommand\thesection{\S\arabic{szn}}

\makeatletter

\renewcommand\l@section{\@dottedtocline{1}{0em}{1.5em}}

\newlength\regolo

\newcommand\sezione[2][]{\stepcounter{szn}\settowidth\regolo{#1}\ifdim\regolo>0pt\@startsection{section}{1}{0pt}{-2\baselineskip}{\baselineskip}{\bfseries\Large}[#1]{#2}\else\@startsection{section}{1}{0pt}{-2\baselineskip}{\baselineskip}{\bfseries\Large}[#2]{#2}\fi\addtocounter{szn}{-1}\refstepcounter{szn}}

\newcommand\appendice[2][]{\refstepcounter{section}\section*{Appendix~\thesection: #2}\settowidth\regolo{#1}\ifdim\regolo>0pt\addcontentsline{toc}{section}{Appendix~\thesection: #1}\else\addcontentsline{toc}{section}{Appendix~\thesection: #2}\fi}

\def\sottosezione#1{\@startsection{subparagraph}{2}{0pt}{-\baselineskip}{.5\baselineskip}{\bfseries}{#1}}

\makeatother

\newenvironment{paragrafetto}[1][]{\refstepcounter{env}\begin{list}{}{\setlength\itemindent{0pt}\setlength\labelsep\parindent\setlength\labelwidth{0pt}\setlength\leftmargin{0pt}\setlength\listparindent\parindent\setlength\parsep\parskip\setlength\partopsep{0pt}}\item\textbf{\arabic{szn}.\arabic{env}~~}\settowidth\regolo{#1}\ifdim\regolo>0pt\textbf{#1\ \ }\else\fi}{\end{list}}

\newenvironment{enunciato}[2][]{\refstepcounter{env}\begin{list}{}{\setlength\itemindent{0pt}\setlength\labelsep{0pt}\setlength\labelwidth{0pt}\setlength\leftmargin\parindent\setlength\listparindent\parindent\setlength\parsep\parskip\setlength\partopsep{0pt}}\item\textbf{\arabic{szn}.\arabic{env}~~}\settowidth\regolo{#2}\ifdim\regolo>0pt\textbf{#2~~}\else\fi\settowidth\regolo{#1}\ifdim\regolo>0pt(#1)~~\else\fi\slshape}{\end{list}}

\newenvironment{enunciato*}[2][]{\begin{list}{}{\setlength\itemindent{0pt}\setlength\labelsep{0pt}\setlength\labelwidth{0pt}\setlength\leftmargin\parindent\setlength\listparindent\parindent\setlength\parsep\parskip\setlength\partopsep{0pt}}\item\settowidth\regolo{#2}\ifdim\regolo>0pt\textbf{#2~~}\else\fi\settowidth\regolo{#1}\ifdim\regolo>0pt(#1)~~\else\fi\slshape}{\end{list}}

\newenvironment{definizione}{\refstepcounter{env}\begin{list}{}{\setlength\itemindent{0pt}\setlength\labelsep{0pt}\setlength\labelwidth{0pt}\setlength\leftmargin{0pt}\setlength\listparindent\parindent\setlength\parsep\parskip\setlength\partopsep{0pt}}\item\textbf{\arabic{szn}.\arabic{env}~~Definition~~}}{\end{list}}

\newenvironment{lemma}[1][]{\begin{enunciato}[#1]{Lemma}}{\end{enunciato}}
\newenvironment{proposizione}[1][]{\begin{enunciato}[#1]{Proposition}}{\end{enunciato}}
\newenvironment{teorema}[1][]{\begin{enunciato}[#1]{Theorem}}{\end{enunciato}}
\newenvironment{corollario}[1][]{\begin{enunciato}[#1]{Corollary}}{\end{enunciato}}

\newenvironment{equazione}[1][]{\equation}{\endequation}
\newenvironment{equazione+}[1][]{\flalign&\phantom|&&}{\endflalign}
\newenvironment{multiriga}[1][]{\multline}{\endmultline}

\providecommand\qedsymbol{\textit{q.e.d.}}
\providecommand\sdsymbol{\ensuremath{\blacksquare}}
\newcommand\mathqed{\quad\hbox{\qedsymbol}}
\newcommand\mathsd{\quad\hbox{\sdsymbol}}

\DeclareRobustCommand\qed{\ifmmode\mathqed\else\leavevmode\unskip\penalty9999\hbox{}\nobreak\hfill\quad\hbox{\qedsymbol}\fi}

\DeclareRobustCommand\sd{\ifmmode\mathsd\else\leavevmode\unskip\penalty9999\hbox{}\nobreak\hfill\quad\hbox{\sdsymbol}\fi}

\newenvironment{proof}[1][]{\begin{list}{}{\setlength\itemindent\parindent\setlength\labelsep{0pt}\setlength\labelwidth{0pt}\setlength\leftmargin{0pt}\setlength\listparindent\parindent\setlength\parsep\parskip\setlength\partopsep{0pt}}\item\settowidth\regolo{#1}\ifdim\regolo>0pt\emph{#1.}\ \else\textit{Proof.~~}\fi}{\qed\end{list}}

\newcounter{itm}

\newenvironment{elenco}{\begin{list}{(\roman{itm})}{\setlength\itemindent{0pt}\setlength\labelsep{0.5em}\setlength\labelwidth\leftmargin\addtolength\labelwidth{-\labelsep}\setlength\listparindent{.5\parindent}\setlength\parsep\parskip\setlength\itemsep\medskipamount\setlength\partopsep{0pt}\usecounter{itm}}}{\end{list}}

\usepackage{amssymb}
\usepackage{mathrsfs}
\usepackage{bbold}

\newcommand\bydef{\stackrel{\text{\tiny\textrm{def}}}{=}}

\newcommand\inciso[1]{\nobreakdash---\hspace{0pt}{#1}\nobreakdash---\hspace{0pt}}

\newcommand\textcap[3][]{\ensuremath{\overset{#1}{\underset{#2}{\textstyle\cap}}\,#3}}

\newcommand\nsubset{\mathrel{\mbox{$\subset$\makebox[0pt][r]{$/\mspace{2mu}$}}}}

\newcommand\nZ{\ensuremath{\mathbb{Z}}}

\newcommand\nR{\ensuremath{\mathbb{R}}}
\newcommand\nC{\ensuremath{\mathbb{C}}}

\newcommand\fro\leftarrow
\newcommand\infro\hookleftarrow
\newcommand\into\hookrightarrow

\newcommand\isoto{\stackrel\thicksim\to}
\newcommand\longfro\longleftarrow
\newcommand\longto\longrightarrow
\newcommand\onfro\twoheadleftarrow
\newcommand\onto\twoheadrightarrow
\newcommand\xfro[1]{\xleftarrow{#1}}
\newcommand\xto[1]{\xrightarrow{#1}}

\newcommand\can\cong
\newcommand\iso\approx

\newcommand\id{\mathit{id}}
\newcommand\pr{\mathit{pr}}

\newcommand\algdim[2][]{\settowidth{\regolo}{\ensuremath{#1}}\ifdim\regolo>0pt\ensuremath{\mathrm{dim}_{#1}#2}\else\ensuremath{\mathrm{dim}\,#2}\fi}
\newcommand\image[1]{\ensuremath{\mathrm{Im}\,#1}}
\newcommand\kernel[1]{\ensuremath{\mathrm{Ker}\,#1}}

\newcommand\Iso[3][]{\ensuremath{\mathrm{Iso}_{#1}(#2,#3)}}
\newcommand\End[2][]{\ensuremath{\mathrm{End}_{#1}(#2)}}
\newcommand\Aut[2][]{\ensuremath{\mathrm{Aut}_{#1}(#2)}}

\newcommand\GL[2][]{\ensuremath{\mathit{GL}_{#1}(#2)}}

\newcommand\supp[2][]{\settowidth{\regolo}{\ensuremath{#1}}\ifdim\regolo>0pt\ensuremath{\mathrm{supp}_{#1}#2}\else\ensuremath{\mathrm{supp}\,#2}\fi}

\newcommand\norm[2][]{\settowidth{\regolo}{\ensuremath{#1}}\ifdim\regolo>0pt\ensuremath{\left\|#2\right\|_{#1}}\else\ensuremath{\left|#2\right|}\fi}

\newcommand\C{\mathit C}
\newcommand\D{\mathit D}

\newcommand\T[2][]{\settowidth{\regolo}{\ensuremath{#1}}\ifdim\regolo>0pt\ensuremath{\mathit T_{#1}#2}\else\ensuremath{\mathit T\,#2}\fi}

\newcommand\aeq\Leftrightarrow
\newcommand\seq\Rightarrow

\newcommand\s[1][]{\settowidth{\regolo}{\ensuremath{#1}}\ifdim\regolo>0pt\ensuremath{{\boldsymbol s}\,#1}\else\ensuremath{\boldsymbol s}\fi}

\renewcommand\t[1][]{\settowidth{\regolo}{\ensuremath{#1}}\ifdim\regolo>0pt\ensuremath{{\boldsymbol t}\,#1}\else\ensuremath{\boldsymbol t}\fi}

\def\ftimes#1#2{\operatorname{_{\mathnormal{#1}}\times_{\mathnormal{#2}}}}

\newcommand\mca[2][]{\settowidth{\regolo}{\ensuremath{#1}}\ifdim\regolo>0pt\ensuremath{#2^{\scriptscriptstyle(#1)}}\else\ensuremath{#2\ftimes{\s}{\t}#2}\fi}

\newcommand\normtang[2][]{\settowidth{\regolo}{\ensuremath{#1}}\ifdim\regolo>0pt\ensuremath{\widetilde{\mathit T}_{#1}#2}\else\ensuremath{\widetilde{\mathit T}\,#2}\fi}

\newcommand\Id{\mathit{Id}}
\newcommand\Ob{\mathrm{Ob}}
\newcommand\Nat{\mathrm{Nat}}

\def\veb#1{\settowidth{\regolo}{\ensuremath{#1}}\ifdim\regolo>0pt\ensuremath{\underline{\mathsf{Vec}}(#1)}\else\ensuremath{\underline{\mathsf{Vec}}}\fi}

\def\rep#1{\ensuremath{\underline{\smash[b]{\mathsf{Rep}}}(#1)}}

\newcommand\forget[1]{\ensuremath{\mathsf F_{#1}}}

\newcommand\bidual[1]{\ensuremath{\mathcal T({#1})}}

\newcommand\eval[1][]{\settowidth{\regolo}{\ensuremath{#1}}\ifdim\regolo>0pt\ensuremath{{\boldsymbol\varepsilon}_{#1}}\else\ensuremath{\boldsymbol\varepsilon}\fi}

\newcommand\R{\ensuremath{\mathscr R}}

\newcommand\canhom[1][]{\settowidth{\regolo}{\ensuremath{#1}}\ifdim\regolo>0pt\ensuremath{{\boldsymbol\pi}_{#1}}\else\ensuremath{\boldsymbol\pi}\fi}

\begin{document}

\title{On the role of\\Effective representations of Lie groupoids\thanks{Research partly supported by Deutsche Forschungs\-gemeinschaft (DFG, German Research Council) in the context of a strategic research proposal by the University of G\"ottingen, and by Fon\-da\-zio\-ne "Ing.~Al\-do Gi\-ni"}}
\author{Giorgio Trentinaglia\thanks{Present affiliation: Courant Research Centre "Higher Order Structures in Mathematics", Georg-August-Universit\"at G\"ottingen, Germany}}
\date{\itshape\footnotesize Department of Mathematics, Utrecht University, The Netherlands}
\maketitle

\begin{abstract}
\noindent In this paper, we undertake the study of the Tannaka duality construction for the ordinary representations of a proper Lie groupoid on vector bundles. We show that for each proper Lie groupoid $\mathcal G$, the canonical homomorphism of $\mathcal G$ into the reconstructed groupoid \bidual{\mathcal G} is surjective, although\inciso{contrary to what happens in the case of groups}it may fail to be an isomorphism. We obtain necessary and sufficient conditions in order that $\mathcal G$ may be isomorphic to \bidual{\mathcal G} and, more generally, in order that \bidual{\mathcal G} may be a Lie groupoid. We show that if \bidual{\mathcal G} is a Lie groupoid, the canonical homomorphism $\mathcal G \to \bidual{\mathcal G}$ is a submersion and the two groupoids have isomorphic categories of representations.
\end{abstract}

\section*{Introduction}
Although Lie groupoids have been intensively studied in the past decades \cite{Co94,Lm98,CdSW,DuZu05}, a satisfactory understanding of their representation theory has not yet been achieved. Notoriously, there are a few fundamental differences between the classical representation theory of Lie groups and its groupoid counterpart, which make it hard to generalize the standard results of the former theory to the latter \cite{2008,AbCr09}. In order to overcome these difficulties, various alternative notions of representation have been proposed recently for Lie groupoids \cite{2008,Abad08}, each of which enables one to attain some specific purpose. Still, in some cases, for example in connection with the study of the presentation problem for ineffective orbifolds \cite{Moe02,HeMe04}, one is by necessity led back into the theory of ordinary representations of Lie groupoids on smooth vector bundles. Very little was known about such representations in general, so far, besides the fact that there are serious difficulties inherent in their construction. (There is no analogue of the \mbox{Peter}--\mbox{Weyl} theorem for them.) In this paper, we begin to remedy this situation by working out the counterpart of the classical duality theory of Tannaka \cite{Tann39,BtD85} in the setting of ordinary representations of Lie groupoids on vector bundles, for proper Lie groupoids.

The classical duality theory of Tannaka addresses the problem of recovering a topological group from its representation ring. In modern terms, this theory can be outlined as follows \cite{JoSt91}. For any such group $G$, the \textit{self-conjugate tensor preserving} natural endomorphisms of the forgetful functor
\begin{equazione}\label{def:Forget}
\mathsf F_G: \rep{G} \xto{\qquad} \veb{}\text, \qquad (V,\rho) \mapsto V
\end{equazione}
form a group $\mathit T(G)$; here \rep{G} stands for the category of continuous, finite dimensional complex $G$-modules and equivariant linear maps between them, and \veb{} stands for the category of finite dimensional complex vector spaces; a natural endomorphism $\lambda \in \End{\mathsf F_G}$ is said to be \textit{tensor preserving} if $\lambda(\nC) =$ identity and $\lambda(R \otimes R') = {\lambda(R) \otimes \lambda(R')}$ for all objects $R$, $R'$ of \rep{G}; \textit{self-conjugate,} if $\lambda\bigl(\overline R\bigr) = \overline{\lambda(R)}$ for all $R$. One endows $\mathit T(G)$ with the smallest topology which makes all the evaluation maps $\lambda \mapsto \lambda(R)$ continuous. There is a continuous group homomorphism \canhom[G] from $G$ into $\mathit T(G)$, defined by setting $\canhom[G](g)(R) = \rho(g)$ for all $R = (V,\rho)$ and all $g \in G$. The fundamental duality theorem of Tannaka states that \canhom[G] is an isomorphism of topological groups whenever the group $G$ is compact Hausdorff.

The same construction can be carried out for a general Lie groupoid $\mathcal G$, if one considers the groupoid analogue of the forgetful functor (\ref{def:Forget}), namely, the functor $\mathsf F_{\mathcal G}$ from the category \rep{\mathcal G} of representations of $\mathcal G$ on vector bundles into the category \veb{M} of complex vector bundles over the base $M$ of $\mathcal G$ which sends $(E,\varrho) \mapsto E$. (The precise details of this construction can be found in Section \ref{sez:BIDUAL} below.) One gets a topological groupoid \bidual{\mathcal G} over $M$, which we shall call the \textit{Tannakian bidual} of the groupoid $\mathcal G$, and a homomorphism of topological groupoids $\canhom[\mathcal G]: \mathcal G \to \bidual{\mathcal G}$ which induces the identity map on the base. It is then natural to ask oneself whether the above-mentioned duality theorem carries over to {\em proper} groupoids (these are the analogue of compact Hausdorff groups in the realm of groupoids). It turns out that this cannot be the case; a standard example, clearly illustrating why, shall be discussed in the second section of this paper.

Despite the lack of a Tannaka duality theorem, we contend that one can still work out, for proper Lie groupoids, a rich and interesting duality theory. In the present paper, we lay down the general framework for this new theory and prove its main theorems. The reasons why we restrict our attention to proper {\em Lie} groupoids are technical, and will become apparent in the course of our exposition.

We proceed to summarize our main contributions:

\smallskip

{\bfseries Theorem \ref{thm:SURJECT}.} {\em For any proper Lie groupoid $\mathcal G$, the canonical homomorphism $\canhom[\mathcal G]: \mathcal G \to \bidual{\mathcal G}$ is full, i.e., surjective.}

\smallskip

\textbullet\ We observe that the Tannakian bidual \bidual{\mathcal G} of any Lie groupoid $\mathcal G$ can be equipped with a natural ``smooth structure'' compatible with its topological groupoid structure, although there is no a priori evidence that this ``smooth structure'' will make \bidual{\mathcal G} into a Lie groupoid. Apart from this, the canonical homomorphism $\canhom[\mathcal G]: \mathcal G \to \bidual{\mathcal G}$ is shown to be always ``smooth''. (The precise definitions are given in Section \ref{sez:BIDUAL}.) The presence of a god-given ``smooth structure'' allows one to assign definite meaning to the question of whether \bidual{\mathcal G} is a Lie groupoid and of whether \canhom[\mathcal G] is an isomorphism of Lie groupoids.

\smallskip

\textbullet\ We say that a Lie groupoid $\mathcal G$ is \textit{reflexive} if the corresponding canonical homomorphism \canhom[\mathcal G] is an isomorphism of ``smooth'' groupoids. Then we prove

\smallskip

{\bfseries Theorem \ref{thm:REFLEX}.} {\em A proper Lie groupoid $\mathcal G$ is reflexive if, and only if, for each base point $x$ there is a representation $\varrho: \mathcal G \to \GL{E}$ such that the kernel of the induced isotropy homomorphism $\varrho_x: \mathcal G_x \to \GL{E_x}$ is trivial.}

\smallskip

\textbullet\ One says that an isotropic arrow $g \in \mathcal G_x$ is \textit{ineffective} if the action of $g$ on the normal space to the $\mathcal G$-orbit at $x$ is trivial. (See Section \ref{sez:MC}.) We have the following characterization of proper Lie groupoids whose bidual is a Lie groupoid (note the formal analogy with the preceding statement):

\smallskip

{\bfseries Theorem \ref{MC}.} {\em For a proper Lie groupoid $\mathcal G$, the following condition is a necessary and also a sufficient one in order that the bidual \bidual{\mathcal G} may be a Lie groupoid: for each base point $x$, there is a representation $\varrho: \mathcal G \to \GL{E}$ such that the kernel of the corresponding isotropy homomorphism $\varrho_x: \mathcal G_x \to \GL{E_x}$ is ineffective, i.e., sits inside the ineffective subgroup of $\mathcal G_x$.}

\smallskip

We shall say, of a representation like the one in the previous statement, that it is \textit{effective at $x$}. We can rephrase the last theorem by saying that the Tannakian bidual of a proper Lie groupoid $\mathcal G$ is a Lie groupoid if, and only if, $\mathcal G$~\textit{has enough effective representations.} The title of the paper refers precisely to this result, which is essentially new to \cite{2008}.

\subsection*{Acknowledgements}

I would like to thank \mbox{I.~Moerdijk} for suggesting the original research problem and for several useful conversations and comments in the early stages of the present work. I am also indebted to \mbox{M.~Crainic} for discussions concerning Haar systems and to \mbox{N.~T.~Zung} for clarifications about his linearizability result.

\tableofcontents

\sezione{The Tannakian bidual of a Lie groupoid}\label{sez:BIDUAL}
We begin this section by introducing some terminological and notational conventions that will be in force throughout the rest of the paper. Some useful references in this connection are the standard books on Lie groupoids \cite{MoeMrc03,Macken05} and the book by \mbox{Bredon} \cite{Bre72}. We then proceed to explain the construction of the Tannakian bidual of a Lie groupoid. In particular, we discuss the associated $\C^\infty$-groupoid structure and prove its basic properties. We conclude with a rapid review of homomorphisms and Morita equivalences.

Roughly speaking, a \textit{Lie groupoid} is an internal groupoid in the category of manifolds of class $\C^\infty$ (recall that the term \textit{groupoid} indicates a small category all morphisms of which are invertible). Thus, a Lie groupoid $\mathcal G$ is given by a pair of such manifolds, \mca[0]{\mathcal G} (the manifold of objects) and \mca[1]{\mathcal G} (the manifold of arrows), a pair of smooth maps $(\s,\t): \mca[1]{\mathcal G} \to {\mca[0]{\mathcal G} \times \mca[0]{\mathcal G}}$ (which are respectively called \textit{source} and \textit{target}) such that the fibred product \mca{\mca[1]{\mathcal G}} exists in the category of $\C^\infty$-manifolds, and a smooth \textit{composition law $\mca{\mca[1]{\mathcal G}} \to \mca[1]{\mathcal G}$} subject to certain requirements. In the first place, there is the requirement that these data should define a category, more precisely, a groupoid. One thereby obtains two uniquely determined maps, the \textit{unit $\mca[0]{\mathcal G} \to \mca[1]{\mathcal G}$} and the \textit{inverse $\mca[1]{\mathcal G} \to \mca[1]{\mathcal G}$}, on which, then, one imposes further smoothness conditions. Moreover \cite{MoeMrc03}, one requires the source map $\s: \mca[1]{\mathcal G} \to \mca[0]{\mathcal G}$ to be a submersion with Hausdorff fibres and the manifold of objects \mca[0]{\mathcal G} to be paracompact (that is, Hausdorff and second countable; cfr.\ \cite{La01}).

Let $\mathcal G = \bigl( \mca[0]{\mathcal G}, \mca[1]{\mathcal G}, \s, \t, \ldots \bigr)$ be a Lie groupoid. The manifold \mca[0]{\mathcal G} is usually called the \textit{base,} and its points, the \textit{base points,} of the groupoid $\mathcal G$. We shall write $\mathcal G$ in place of \mca[1]{\mathcal G}, as this is not likely to cause any confusion. For each arrow $g \in \mathcal G$, the base points $\s(g)$ and $\t(g)$ are resp.\ called the \textit{source} and the \textit{target} of $g$. We shall not distinguish notationally between a base point and the corresponding unit arrow. We shall adopt the abbreviations
\begin{equazione}\label{not:G(S,S')}
\mathcal G(S,S') := \bigl\{g \in \mathcal G: \s(g) \in S \text{~\&~} \t(g) \in S'\bigr\}\text, \qquad \mathcal G_S := \mathcal G(S,S)\text,
\end{equazione}
and $\mathcal G^S := \mathcal G(S,\text-) := \mathcal G\bigl(S,\mca[0]{\mathcal G}\bigr) = \s^{-1}(S)$ for all subsets $S, S' \subset \mca[0]{\mathcal G}$. When $S = \{x\}$ and $S' = \{x'\}$ are singletons, we further abbreviate the respective instances of (\ref{not:G(S,S')}) into $\mathcal G(x,x')$, $\mathcal G_x$ and $\mathcal G^x$. One refers to $\mathcal G_x$ as the \textit{isotropy group of~$\,\mathcal G$ at $x$}. This is in fact a Lie group, embedded into $\mathcal G$ as a closed submanifold.

\begin{paragrafetto}[Fundamental example]\label{xmp:GL(E)}
Any vector bundle $E$ of class $\C^\infty$ (real or complex, of globally constant rank) over a smooth manifold $M$ determines a Lie groupoid \GL{E} with base $M$, called the \textit{linear groupoid associated with $E$}. The arrows $\in \GL{E}(x,x')$ are the linear isomorphisms $E_x \isoto E_{x'}$ between the corresponding fibres of the vector bundle $E$.
\end{paragrafetto}

A \textit{homomorphism of Lie groupoids $\phi: \mathcal G \to \mathcal H$} is defined to be a smooth functor $\mca[0]{\phi}: \mca[0]{\mathcal G} \to \mca[0]{\mathcal H}$, $\mca[1]{\phi}: \mca[1]{\mathcal G} \to \mca[1]{\mathcal H}$.

\begin{definizione}\label{def:REPR}
Let $\mathcal G$ be a Lie groupoid, and let $E$ be a vector bundle over \mca[0]{\mathcal G}. By a \textit{representation of~$\,\mathcal G$ on $E$}, we mean a Lie groupoid homomorphism $\varrho: \mathcal G \to \GL{E}$ that induces the identity on the base. Thus, a representation $\varrho$ assigns each arrow $g \in \mathcal G(x,x')$ a linear isomorphism $\varrho(g): E_x \isoto E_{x'}$ in a smooth and functorial way. More pedantically, we shall think of representations as pairs $(E,\varrho)$ formed by a vector bundle $E$ and a homomorphism $\varrho$ as above.
\end{definizione}

We regard the notions to be discussed next in \S\ref{C^infty-Sp} and \S\ref{C^infty-Gpd} as standard. No attempts to make comparisons with the literature shall be made here, nor claims to originality. In any case, we are not interested in developing a general theory. This said, we are entitled to use these notions freely for our purposes.

\begin{paragrafetto}[$\C^\infty$-Spaces]\label{C^infty-Sp}
Recall that a \textit{functionally structured space} is a topological space $X$ endowed with a sheaf $\mathscr F_X$ of real algebras of continuous real valued functions on $X$ (functional structure). A \textit{smooth} mapping from a functionally structured space $(X,\mathscr F_X)$ into another such space $(Y,\mathscr F_Y)$ is a continuous mapping $f: X \to Y$ such that ${\alpha \circ f} \in \mathscr F_X(f^{-1}(V))$ for every open subset $V \subset Y$ and every function $\alpha \in \mathscr F_Y(V)$. (Compare \cite{Bre72}, p.~297.)

Let $\mathscr F$ be a given functional structure on a topological space $X$. We shall let $\mathscr F^\infty$ denote the sheaf of continuous real valued functions on $X$ generated by the following presheaf (of such functions)
\begin{multiriga}\label{C^infty-Sp1}
U\: \mapsto\: \bigl\{\, f(a_1|_U, \ldots, a_d|_U)\: \bigl|\: \text{~$f: \nR^d \to \nR$ of class $\C^\infty$,}\\
a_1, \ldots, a_d \in \mathscr F(U)\, \bigr\}
\end{multiriga}
where the expression $f(a_1|_U,\ldots,a_d|_U)$ indicates, of course, the function $u \mapsto f(a_1(u),\ldots,a_d(u))$ on $U$. By a \textit{$\C^\infty$-space,} we mean a functionally structured space $X$ such that $\mathscr F_X = {\mathscr F_X}^\infty$. We then say that $\mathscr F_X$ is a \textit{$\C^\infty$-functional structure} on $X$. Note that smooth manifolds can be regarded as topological spaces endowed with a $\C^\infty$-functional structure locally isomorphic to that given by the smooth functions on $\nR^n$. $\C^\infty$-Spaces and their smooth mappings (\textit{$\C^\infty$-mappings,} for short) form a category which, for certain purposes, is more convenient than the full subcategory formed by smooth manifolds.
\end{paragrafetto}

\begin{paragrafetto}[$\C^\infty$-Groupoids]\label{C^infty-Gpd}
Observe that if $(X,\mathscr F)$ is a $\C^\infty$-space, so is $(S,\mathscr F|_S)$ for any subspace $S$ of $X$; here $\mathscr F|_S := {{i_S}^*\mathscr F}$ denotes the functional structure on $S$ induced by $\mathscr F$ along the inclusion $i_S: S \into X$. (Recall that for an arbitrary continuous mapping $f: S \to T$ into a functionally structured space $(T,\mathscr T)$, ${f^*\mathscr T}$ denotes the functional sheaf on $S$ formed by the functions which are locally the pullback along $f$ of functions in $\mathscr T$. Compare \cite{Bre72}, p.~297 again.) The induced $\C^\infty$-structure $(S,\mathscr F|_S)$ has the following {\em universal property:} for any commutative diagram of maps between $\C^\infty$-spaces
\begin{equazione}\label{UNIVPROP:S->X}\begin{split}
\xymatrix@C=80pt@M=5pt{(Y,\mathscr G) \ar@{-->}[dr]_{f'} \ar[r]^-f & (X,\mathscr F) \\ & (S,\mathscr F|_S)\text,\!\! \ar@{^(->}[u]_{i_S}}
\end{split}\end{equazione}
the map $f$ is smooth if, and only if, the same is true of the map $f'$.

Observe, next, that if $(X,\mathscr F)$ and $(Y,\mathscr G)$ are arbitrary functionally structured spaces, then so is their Cartesian product endowed with the sheaf ${\mathscr F \otimes \mathscr G}$ locally generated by the functions $({\varphi \otimes \psi})(x,y) = \varphi(x)\psi(y)$. It follows that $({\mathscr F^\infty \otimes \mathscr G^\infty})^\infty$ is a $\C^\infty$-functional structure on $X\times Y$, turning this into the product of $(X,\mathscr F^\infty)$ and $(Y,\mathscr G^\infty)$ in the category of $\C^\infty$-spaces. Taking the universal property (\ref{UNIVPROP:S->X}) into account, one sees that the category of $\C^\infty$-spaces is closed under fibred products (pullbacks). Notice that when $X$ and $Y$ are smooth manifolds, or when $S \subset X$ is a submanifold, one recovers the correct manifold structures, so that all these constructions for $\C^\infty$-spaces agree with the usual ones on manifolds whenever the latter make sense.

We shall use the term \textit{$\C^\infty$-groupoid} to indicate a groupoid whose set of objects and of arrows are each endowed with the structure of a $\C^\infty$-space so that all the maps arising from the groupoid structure (source, target, composition, unit section, inverse) are morphisms of $\C^\infty$-spaces. The base space $X$ will always be a smooth manifold in practice, with $\C^\infty$ functional structure given by the sheaf of smooth functions on $X$. Every Lie groupoid is, in particular, an example of a $\C^\infty$-groupoid.
\end{paragrafetto}

\sottosezione{Tannakian biduals and representative functions}

The complex linear representations of a Lie groupoid $\mathcal G$ shall be regarded as the objects of a category, which shall be denoted by \rep{\mathcal G} hereafter. The morphisms $(E,\varrho) \to (E',\varrho')$ in \rep{\mathcal G} are given by the \textit{intertwiners,} i.e., those morphisms $a: E \to E'$ of complex linear vector bundles over \mca[0]{\mathcal G} which are compatible with the given actions of $\mathcal G$ in the sense that ${a_{x'} \circ \varrho(g)} = {\varrho'(g) \circ a_x}$ $\forall g \in \mathcal G(x,x')$. In the case of groups, one recovers the usual equivariant linear maps. The familiar operations of Representation Theory\inciso{direct sum, tensor product, complex conjugation, contragredient (or dual) etc.}lend themselves to an obvious generalization; one uses the corresponding standard operations on complex linear vector bundles.

Let $M$ be the base of $\mathcal G$, and let $x_0 \in M$. Consider the functor
\begin{equazione}\label{N.v1}
{x_0}^*: \veb{M} \longto \veb{}\text, \qquad E \mapsto E_{x_0}
\end{equazione}
which assigns each vector bundle $E$ over $M$ its fibre at the point $x_0$, and let \forget{\mathcal G,x_0} denote the composite
\begin{equazione}\label{N.v2}
\rep{\mathcal G} \xto{\; \forget{\mathcal G} \;} \veb{M} \xto{\; {x_0}^* \;} \veb{}\text, \qquad (E,\varrho) \mapsto E \mapsto E_{x_0}\text.
\end{equazione}

\begin{definizione}\label{def:BID}
The \textit{Tannakian bidual} of $\mathcal G$ is the groupoid \bidual{\mathcal G} over $M$ defined as follows. For each pair of base points $x, x' \in M$, put
\begin{equazione}\label{N.v3}
\bidual{\mathcal G}(x,x') := \Nat^{\overline\otimes}\bigl(\forget{\mathcal G,x},\forget{\mathcal G,x'}\bigr)\text;
\end{equazione}
the right-hand side here denotes the set of all self-conjugate, tensor preserving natural transformations from $\forget{\mathcal G,x}$ to $\forget{\mathcal G,x'}$, that is, natural transformations $\lambda: \forget{\mathcal G,x} \to \forget{\mathcal G,x'}$ such that the following diagrams commute
\begin{equazione}\label{def:TP,SC}
\xymatrix@C=25pt@R=20pt{{E_x \otimes F_x} \ar@{<->}[d]^\simeq \ar[rr]^-{\lambda(R)\,\otimes\,\lambda(S)} & & {E_{x'} \otimes F_{x'}} \ar@{<->}[d]^\simeq & \nC \ar@{<->}[d]^\simeq \ar@{=}[r] & \nC \ar@{<->}[d]^\simeq & \overline{(E_x)} \ar@{<->}[d]^\simeq \ar[r]^-{\overline{\lambda(R)}} & \overline{(E_{x'})} \ar@{<->}[d]^\simeq \\ (E\otimes F)_x \ar[rr]^-{\lambda(R\,\otimes\,S)} & & (E\otimes F)_{x'} & \underline \nC_x \ar[r]^-{\lambda(\underline\nC)} & \underline \nC_{x'}\negthickspace & \overline E_x \ar[r]^-{\lambda(\overline R)} & \overline E_{x'}\negthickspace}
\end{equazione}
for all $R = (E,\varrho)$, $S = (F,\varsigma) \in {\Ob\:\rep{\mathcal G}}$. As to the groupoid structure, the composition law is defined to be $({\lambda' \cdot \lambda})(R) := {\lambda'(R) \circ \lambda(R)}$.
\end{definizione}

\begin{paragrafetto}[Remarks]\label{small,Nat=Iso}
\bidual{\mathcal G} is a small category, because, clearly, the category \rep{\mathcal G} possesses a small skeleton; compare \cite{MacLane71}, p.~91.

From the rigidity of the tensor category \rep{\mathcal G}\inciso{that is to say, roughly speaking, from the existence of duals in \rep{\mathcal G}}it follows that any tensor preserving natural transformation $\forget{\mathcal G,x} \to \forget{\mathcal G,x'}$ is necessarily an {\em isomorphism.} Hence \bidual{\mathcal G} is really a {\em groupoid.} For the necessary explanations and a proof of this fact, we refer the reader to \cite{DeMi82}, p.~117.
\end{paragrafetto}

An object $R = (E,\varrho)$ of the category \rep{\mathcal G} determines a homomorphism of groupoids over $M$ (identical on the base manifold $M$)
\begin{equazione}\label{def:EV_R}
\eval[R]: \bidual{\mathcal G} \longto \GL{E}, \qquad \lambda \mapsto \lambda(R)\text.
\end{equazione}
We shall call \eval[R] the \textit{evaluation representation at $R$}.

Let $\phi$ be any Hilbert metric on the (smooth, complex) vector bundle $E$. Let $\zeta$ and $\zeta'$ be any smooth sections of $E$ over $M$. Consider the following function on (the manifold of arrows of) the linear groupoid \GL{E}:
\begin{equazione}\label{N.v8}
q_{\phi,\zeta,\zeta'}: \GL E \to \nC\text, \qquad \GL E(x,x') \ni \mu \mapsto \bigl\langle {\mu \cdot \zeta(x)}, {\zeta'(x')} \bigr\rangle\text.
\end{equazione}
Define $\R$ as the collection of all the functions on (the set of arrows of) the Tannakian bidual \bidual{\mathcal G} which can be written as
\begin{equazione}\label{N.v9}
r_{R,\phi,\zeta,\zeta'} := {q_{\phi,\zeta,\zeta'} \circ \eval[R]}
\end{equazione}
for some $R = (E,\varrho) \in {\Ob\:\rep{\mathcal G}}$, some Hilbert metric $\phi$ on $E$, and some global sections $\zeta$, $\zeta'$ as above. We call the elements of $\R$ \textit{representative functions.} It is easy to see that $\R$ is a complex algebra of functions on the set of arrows of the groupoid \bidual{\mathcal G}, closed under complex conjugation. It follows that the real and imaginary parts of any function in $\R$ also belong to $\R$.

We endow the set of arrows of the groupoid \bidual{\mathcal G} with the smallest topology making all the representative functions continuous. Since Hilbert metrics exist on any smooth vector bundle over a paracompact base manifold, one obtains a Hausdorff topological space. The real valued representative functions generate, on this topological space, a functional structure, which we shall complete to a $\C^\infty$-space structure $\R^\infty$ as explained in \ref{C^infty-Sp}. This $\C^\infty$-space structure can be given the following ``universal'' characterization:

\begin{lemma}\label{lem:UNIVPROP}
There exists a unique $\C^\infty$-space structure on (the set of arrows of) \bidual{\mathcal G} with the following property: for each map $f: X \to \bidual{\mathcal G}$ from a $\C^\infty$-space $(X,\mathscr F_X)$, $f$ is a smooth mapping of $\C^\infty$-spaces if, and only if, the composition
\begin{equazione+}\label{equ:UNIVPROP}
{\eval[R] \circ f}: X \to \GL{E} \quad \text{is smooth} & \forall R = (E,\varrho) \in \Ob\:\rep{\mathcal G}\text.
\end{equazione+}
This $\C^\infty$-space structure is precisely the one defined by the representative functions.
\end{lemma}\begin{proof}
We limit ourselves to checking that the $\C^\infty$-space structure defined by the representative functions has the indicated universal property. So, let a map $f: X \to \bidual{\mathcal G}$ be given.

In one direction, we have to show that the evaluation representation \eval[R] is a smooth mapping, for each $R = (E,\varrho)$. It will be enough to note that if $\phi$ is any Hilbert metric on $E$, then the local functions $q_{\phi,\zeta_i,\zeta_j}$ provide, when $\{\zeta_i\}$ is made to vary over all local frames of sections of the vector bundle $E$, local coordinate systems on the manifold \GL{E} in terms of which any smooth function on the same manifold can be expressed locally.

Conversely, let $r \equiv r_{R,\phi,\zeta,\zeta'}$ be a given representative function; we have to show that (\ref{equ:UNIVPROP}) implies ${r \circ f} \in \mathscr F_X(X)$. This is clear, because ${r \circ f} = {q_{\phi,\zeta,\zeta'} \circ \eval[R] \circ f}$ and $q_{\phi,\zeta,\zeta'}$ is smooth.
\end{proof}

\begin{proposizione}\label{prop:T(G)=CinftyGPD}
With the $\C^\infty$-structure described above, the Tannakian bidual \bidual{\mathcal G} of any Lie groupoid $\mathcal G$ is a $\C^\infty$-groupoid.
\end{proposizione}\begin{proof}
By the preceding lemma, the smoothness of the composition law $\mca{\bidual{\mathcal G}} \to \bidual{\mathcal G}$ and of the inverse $\bidual{\mathcal G} \to \bidual{\mathcal G}$ is a consequence of the commutativity of the following diagrams
$$%
\xymatrix@C=40pt{{\mca{\bidual{\mathcal G}}}\ar[d]^{\eval[R] \times\, \eval[R]}\ar[r] & \bidual{\mathcal G}\ar[d]^{\eval[R]} & \bidual{\mathcal G}\ar[d]^{\eval[R]}\ar[r] & \bidual{\mathcal G}\ar[d]^{\eval[R]} \\ {\mca{\GL E}}\ar[r]^(.60){\smash{\text{compos.}}} & {\GL E} & {\GL E}\ar[r]^(.48){\text{inverse}} & {\GL E}}
$$%
for all $R = (E,\varrho) \in \Ob\:\rep{\mathcal G}$. The smoothness of the other structure maps should be clear.
\end{proof}

Note that there is an obvious notion of \textit{smooth representation of a $\C^\infty$-groupoid on a smooth vector bundle.} Thus, we can make sense of the notation $\underline{\smash{\mathsf{Rep}}}\bigl(\bidual{\mathcal G}\bigr)$ for any Lie groupoid $\mathcal G$. Evaluation (\ref{def:EV_R}) defines a functor
\begin{equazione}\label{def:EVALFUNC}
\eval: \rep{\mathcal G} \longto \underline{\smash{\mathsf{Rep}}}\bigl(\bidual{\mathcal G}\bigr), \qquad R = (E,\varrho) \mapsto (E,\eval[R])\text.
\end{equazione}
It is natural to ask oneself whether this functor is an equivalence of categories; between Sections \ref{sez:CANHOM} and \ref{sez:REPCHARTS}, we shall see that the answer is affirmative whenever $\mathcal G$ is {\em proper,} even though \bidual{\mathcal G} need not be isomorphic to $\mathcal G$.

\sottosezione{Invariance under Morita equivalences}

Given a homomorphism of Lie groupoids $\phi: \mathcal G \to \mathcal H$, we define the \textit{inverse image functor $\phi^*: \rep{\mathcal H} \longto \rep{\mathcal G}$} (which we shall also call \textit{pullback along $\phi$}) as follows. Let $f: M \to N$ denote the map induced by $\phi$ on the bases of the two groupoids.

Let $S = (F,\varsigma)$ be a representation of the groupoid $\mathcal H$. The linear groupoid \GL{f^*F} is the groupoid over $M$ induced by the groupoid $\GL{F} \rightrightarrows N$ along the base change $f: M \to N$. (See \cite{MoeMrc03}.) To put it differently,
$$%
\GL{f^*F} = (f \times f)^*\bigl[\GL{F}\bigr]
$$%
in the category of smooth manifolds. By the universal property of induced groupoids, there is a unique Lie groupoid homomorphism $\phi^*(\varsigma): \mathcal G \to \GL{f^*F}$ such that when one composes it with the canonical homomorphism
\begin{equazione}\label{GL(f*F)->GL(F)}
\gamma: \GL{f^*F} \longto \GL{F}
\end{equazione}
one gets ${\varsigma\circ\phi}: \mathcal G \to \GL{F}$. We define $\phi^*(S) := \bigl(f^*F,\phi^*(\varsigma)\bigr)$. Note that if $b: S \to S' = (F',\varsigma')$ is a morphism in the category \rep{\mathcal H} then the morphism of vector bundles $f^*b: f^*F \to f^*(F')$ intertwines the two $\mathcal G$-actions $\phi^*(\varsigma)$ and $\phi^*(\varsigma')$ and hence defines a morphism $\phi^*(S) \to \phi^*(S')$ in the category \rep{\mathcal G}. This concludes the description of the functor $\phi^*$.

\begin{paragrafetto}[Remark]\label{rmk:T(phi)}
By taking into account the identity ${\forget{\mathcal G} \circ \phi^*} = {f^* \circ \forget{\mathcal H}}$, it is straightforward to see that the inverse image functor $\phi^*$ determines a homomorphism of $\C^\infty$-groupoids $\bidual{\phi}: \bidual{\mathcal H} \to \bidual{\mathcal G}$. (Compare \cite{2008}, \S24.) In fact, the correspondence $(\text-) \mapsto \bidual{\text-}$ yields a {\em functor} from the category of Lie groupoids into that of $\C^\infty$-groupoids.
\end{paragrafetto}

Recall that a homomorphism of Lie groupoids $\phi: \mathcal G \to \mathcal H$ (notations as above) is said to be a \textit{Morita equivalence} if the diagram
\begin{equazione}\label{def:ME1}\begin{split}
\xymatrix@C=30pt{\mca[1]{\mathcal G}\ar[d]_{(\s,\t)}\ar[r]^{\smash{\mca[1]\phi}} & \mca[1]{\mathcal H}\ar[d]^{(\s,\t)} \\ {M\times M}\ar[r]^-{f\times f} & {N\times N}}
\end{split}\end{equazione}
is a pullback in the category of manifolds of class $\C^\infty$ and the mapping
\begin{equazione}\label{def:ME2}
{\mca[1]{\mathcal H} \ftimes{\s}{f} M} \longto N\text, \qquad (h,x) \mapsto \t(h)
\end{equazione}
is a surjective submersion. It is a well known fact that the pullback functor $\phi^*: \rep{\mathcal H} \longto \rep{\mathcal G}$ associated with a Morita equivalence $\phi: \mathcal G \to \mathcal H$ is an equivalence of categories. For the benefit of the reader and for later reference, we shall now briefly review the proof of this result.

We need a preliminary remark. Let $\phi, \psi: \mathcal G \to \mathcal H$ be two arbitrary homomorphisms of Lie groupoids. Recall that a \textit{transformation $\theta: \phi \to \psi$} is an ordinary natural transformation of functors $\phi \to \psi$ with the additional property that, when regarded as a map $\theta: \mca[0]{\mathcal G} \to \mca[1]{\mathcal H}$, it is smooth. For each representation $S = (F,\varsigma)$ of~$\,\mathcal H$, one obtains a smooth section $\theta^*(S)$ of the vector bundle \Iso[M]{\phi^*F}{\psi^*F} over $M \equiv \mca[0]{\mathcal G}$ (in other words, an isomorphism $\theta^*(S): \phi^*F \simeq \psi^*F$ of vector bundles over $M$) by exploiting the universal property of that bundle as the pullback to $M$ of the submersion $\GL{F} \to {\mca[0]{\mathcal H} \times \mca[0]{\mathcal H}}$ along the smooth mapping $x \mapsto \bigl(\phi(x),\psi(x)\bigr)$. The rule $S \mapsto \theta^*(S)$ defines, in fact, a natural isomorphism $\theta^*: \phi^* \simeq \psi^*$.

Let us go back to our Morita equivalence $\phi: \mathcal G \to \mathcal H$. In order to construct a quasi-inverse for the functor $\phi^*$, let us say $\phi_!$, it is not restrictive to assume that the base map $f: M \to N$ associated with $\phi$ is a {\em surjective submersion.} Indeed, if we take the following \textit{weak pullback} (see \cite{MoeMrc03}, pp.~123--132)
$$%
\xymatrix@C=40pt{\mathcal P\ar[d]^\chi\ar[r]^\psi & \mathcal G\ar[d]^\phi_{}="b" \\ \mathcal H\ar@{=}[r]^(.4)\quad="a"\ar@{=>}@/^.5pc/"a";"b"^\theta & \mathcal H}
$$%
then the Lie groupoid homomorphisms $\psi$ and $\chi$ are Morita equivalences with the property that the respective base maps are surjective submersions. Now, if we prove that $\psi^*$ and $\chi^*$ are categorical equivalences then the same will be true of $\phi^*$, because, by the remarks in the previous paragraph, we have natural equivalences $\chi^* \simeq ({\phi \circ \psi})^* \simeq {\psi^* \circ \phi^*}$.

Under the preceding simplifying assumption, there will be an open cover of the manifold $N$ by open subsets $U_i$ such that one can find smooth sections $\alpha_i: U_i \into M$ to the mapping $f$.

Let a representation $R = (E,\varrho) \in {\mathrm{Ob}\: \rep{\mathcal G}}$ be given. One constructs a smooth vector bundle $\phi_!(E)$ in \veb{N}, as follows. Put $E_i := {\alpha_i}^*E \in {\mathrm{Ob}\: \veb{U_i}}$. Introduce the cocycle of vector bundle isomorphisms
$$%
c_{ij}: E_j|U_{ij} \isoto E_i|U_{ij}
$$%
given, for each $x \in U_{ij} := {U_i \cap U_j}$, by the composite linear isomorphism $({\alpha_j}^* E)_x \simeq E_{\alpha_j(x)} \xto{\; \varrho(g) \;} E_{\alpha_i(x)} \simeq ({\alpha_i}^* E)_x$, where $g: \alpha_j(x) \to \alpha_i(x)$ is the unique arrow in $\mathcal G$ such that $\phi(g) = x$. From these data one can construct the required vector bundle $\phi_!(E)$; the elements of its total space are the equivalence classes of triples $(i,x,e)$, with $x \in U_i$ and $e \in E_{\alpha_i(x)}$, the equivalence being
$$%
(i,x,e) \thicksim (j,y,f) \qquad \aeq \qquad x = y \in U_{ij} \quad \text{and} \quad c_{ij}(f) = e\text;
$$%
the projection onto $N$ is given by $[i,x,e] \mapsto x$.

An action $\phi_!(\varrho)$ of the Lie groupoid $\mathcal H$ on this new vector bundle $\phi_!(E)$ can be obtained as follows. Let an arrow $h: x \to x'$ in $\mathcal H$ be given. Choose indices $i$ and $i'$ such that $x \in U_i$ and $x' \in U_{i'}$. There exists a unique arrow $g: \alpha_i(x) \to \alpha_{i'}(x')$ such that $\phi(g) = h$. Then, let
$$%
\phi_!(\varrho)(g): [i,x,e] \mapsto \bigl[i',x',{\varrho(k) \cdot e}\bigr]\text.
$$%

The pair $\phi_!(R) := \bigl(\phi_!(E),\phi_!(\varrho)\bigr)$ is an object of the category \rep{\mathcal G}. We shall leave the rest of the construction as an exercise.

\sezione{The canonical homomorphism}\label{sez:CANHOM}
Recall that a Lie (topological, $\C^\infty$-) groupoid $\mathcal G$ is said to be \textit{proper,} when it is Hausdorff, and the combined source--target map
\begin{equazione}\label{def:(source,target)}
(\s,\t): \mca[1]{\mathcal G} \to {\mca[0]{\mathcal G} \times \mca[0]{\mathcal G}}\text, \qquad g \mapsto (\s(g),\t(g))
\end{equazione}
is proper (in the familiar sense: the inverse image of a compact subset is compact). Compare \cite{MoeMrc03}. When $\mathcal G$ is proper, every $\mathcal G$-orbit is a closed submanifold.

Normalized Haar systems on proper groupoids are the analogue of Haar probability measures on compact groups. The precise definition is as follows:

\begin{definizione}\label{def:NHS}
A \textit{normalized Haar system} (on a Lie groupoid $\mathcal G$ over $M$) is a family $\mu = \{\mu^x: x \in M\}$ of positive Radon measures, each one with support in the respective source fibre $\mathcal G^x$, such that the following conditions are satisfied:
\begin{elenco}
\item all smooth functions on $\mathcal G^x$ are integrable with respect to $\mu^x$, that is to say, $\C^\infty(\mathcal G^x) \subset \mathit L^1(\mu^x)$;
\item(smoothness) for each $\varphi \in \C^\infty(\mca[1]{\mathcal G})$, the function $\Phi$ on $M$ given by
$$%
x \mapsto \Phi(x) \bydef {\int_{\mathcal G^x} \varphi|_{\mathcal G^x}\, \mathnormal d\mu^x} \qquad \text{is of class $\C^\infty$;}
$$%
\item(right invariance) for each $g \in \mathcal G(x,y)$, and for all $\varphi \in \C^\infty(\mathcal G^x)$, one has
$$%
{\int_{\mathcal G^y} {\varphi \circ \tau^g}\, \mathit d\mu^y} = {\int_{\mathcal G^x} \varphi\, \mathit d\mu^x}\text,
$$%
where $\tau^g: \mathcal G(y,\text-) \to \mathcal G(x,\text-)$ denotes right translation by $g$;
\item(normalization) ${\int \mathit d\mu^x} = \mu^x(\mathcal G^x) = 1$, for every $x \in M$.
\end{elenco}
\end{definizione}

It can be shown that every {\em proper} Lie groupoid admits normalized Haar systems. Compare \cite{CraMoe01,Cra03}.

Observe that if $E$ is a smooth vector bundle over the base of the groupoid $\mathcal G$, and $\psi: \mathcal G \to E$ is a smooth mapping such that the source fibre $\mathcal G(x,\text-)$ is mapped into the vector space $E_x$ for every base point $x$, then the integral
$$%
\Psi(x) \bydef {\int \psi|_{\mathcal G^x}\, \mathit d\mu^x}
$$%
defines a global smooth section $\Psi$ to $E$. This follows easily from the previous definition, by working in local coordinates.

\medskip

One can think of any Lie groupoid $\mathcal G$ as acting on its own base manifold $M$. A subset $S \subset M$ will be called \textit{invariant} if ${\mathcal G\cdot s} \subset S$ for all $s \in S$. Let $(E,\varrho)$ be a representation of $\mathcal G$. A partial section $\varphi: S \to E$, defined over an invariant subset $S \subset M$, will be said to be \textit{equivariant with respect to $\varrho$} if $\varphi(\t[g]) = {\varrho(g) \cdot \varphi(\s[g])}$ for all $g \in \mathcal G(S,S)$. We shall call an arbitrary partial section $\varphi: S \to E$ \textit{smooth,} when for every $m \in M$ one can find an open neighbourhood $B \subset M$ of $m$ and a smooth local section to $E$ over $B$ which restricts to $\varphi$ on the intersection ${B \cap S}$.

\begin{enunciato}{Equivariant Extension Lemma}\label{lem:EQVEXT}
Let $\mathcal G$ be a proper Lie groupoid, and let $M$ denote its base manifold. Let $(E,\varrho)$ be a representation of $\mathcal G$. Then, each $\varrho$-equivariant smooth partial section $\varphi: S \to E$ (defined over an invariant subset $S$ of $M$) can be extended to a global $\varrho$-equivariant smooth section $\Phi: M \to E$.
\end{enunciato}\begin{proof}
To begin with, we construct a smooth section $\psi: M \to E$ extending $\varphi$ (possibly not equivariant), as follows. Cover $M$ with a family of open subsets $\{B_i\}$ such that for each $i$ there exists a smooth local section $\varphi_i: B_i \to E$ extending $\varphi$. Since $M$ is paracompact, there is some smooth partition of unity $\{f_j\}$ over $M$, subordinated to the given open cover. Then,
$$%
\psi \bydef {\sum f_j \varphi_{i(j)}}
$$%
is the desired global extension.

The correspondence $g \mapsto {\varrho(g^{-1}) \cdot \psi(\t[g])}$ yields a smooth map $\Psi: \mathcal G \to E$, which, for every $x \in M$, restricts to a smooth map $\Psi^x: \mathcal G^x \to E_x$ on the corresponding source fibre. Since $\mathcal G$ is proper, we can fix a normalized Haar system $\{\mu^x\}$ on $\mathcal G$. For all $x$, the vector valued smooth map $\Psi^x$ must be integrable with respect to $\mu^x$. Then, put
$$%
\Phi(x) \bydef {\int \Psi^x\, \mathrm d\mu^x} \in E_x\text.
$$%
It can be easily checked that the resulting smooth section $\Phi: M \to E$ extends $\varphi$ and is equivariant with respect to $\varrho$.
\end{proof}

\begin{proposizione}\label{prop:EQVEXT}
Let $(E,\varrho)$ and $(F,\varsigma)$ be representations of a proper Lie groupoid $\mathcal G$. Let $x_0$ be a point of the base manifold $M$ of $\mathcal G$, and let $G$ denote the isotropy group of $\mathcal G$ at $x_0$. Suppose that $A$ is a $G$-equivariant linear map of $E_{x_0}$ into $F_{x_0}$. Then there exists a morphism of representations $a: (E,\varrho) \to (F,\varsigma)$ in \rep{\mathcal G} such that $a_{x_0} = A$.
\end{proposizione}\begin{proof}
Let $\mathit L(E,F)$ be the vector bundle over $M$ whose fibre at a generic base point $x$ is the vector space $\mathit L(E_x,F_x)$. Given any arrow $g \in \mathcal G(x,x')$, we let $\mathit L(\varrho,\varsigma)(g)$ denote the linear map
$$%
\bigl\{\lambda \mapsto {\varsigma(g) \circ \lambda \circ \varrho(g^{-1})}\bigr\}: \mathit L(E_x,F_x) \longto \mathit L(E_{x'},F_{x'})\text.
$$%
The pair $\bigl(\mathit L(E,F),\mathit L(\varrho,\varsigma)\bigr)$ is a representation of the groupoid $\mathcal G$. Let $S \subset M$ denote the $\mathcal G$-orbit through $x_0$, and let $\varphi: S \to \mathit L(E,F)$ be the map
$$%
\bigl\{ s\; \mapsto\; {\varsigma(g) \circ A \circ \varrho(g^{-1})}\bigl| g \in \mathcal G(x_0,s) \bigr\}\text.
$$%
(This is well-defined, because $A$ is a homomorphism of $G$-modules.) Note that $\varphi$ is a section of $\mathit L(E,F)$ over $S$, equivariant with respect to $\mathit L(\varrho,\varsigma)$. Furthermore, $\varphi$ is a smooth mapping of $S$ into $\mathit L(E,F)$, because, locally about any point $s_0 \in S$, it is obtained by composing a smooth target section $W_0 \into \mathcal G(x_0,\text-)$ (defined over a small open neighbourhood $W_0 \subset S$ of $s_0$) with the smooth map $\mathcal G(x_0,\text-) \to \mathit L(E,F)$ given by $g \mapsto {\varsigma(g) \circ A \circ \varrho(g^{-1})}$. Since $S$ is a submanifold of $M$, it follows that $\varphi$ is a smooth partial section of $\mathit L(E,F)$.

The Equivariant Extension Lemma then provides us with a global smooth section $a: M \to \mathit L(E,F)$, equivariant with respect to $\mathit L(\varrho,\varsigma)$, and extending $\varphi$. The proof is now finished.
\end{proof}

\begin{definizione}\label{def:CANHOM}
One has a \textit{canonical homomorphism $\canhom[\mathcal G]: \mathcal G \longto \bidual{\mathcal G}$}, defined by setting
\begin{equazione+}\label{equ:CANHOM}
{\canhom[\mathcal G](g)}(R) := \varrho(g) & \text{for all~$\:R = (E,\varrho) \in \Ob\:\rep{\mathcal G}$.}
\end{equazione+}\end{definizione}

From Lemma \ref{lem:UNIVPROP} and the identities ${\eval[(E,\varrho)] \circ \canhom[\mathcal G]} = \varrho$, it follows that \canhom[\mathcal G] is a homomorphism of $\C^\infty$-groupoids. In fact, the canonical homomorphisms altogether determine a natural transformation $\canhom: \Id \to \bidual{\text-}$ (between functors from the category of Lie groupoids into that of $\C^\infty$-groupoids).

\medskip

The next theorem is the central result of this section. For the benefit of the reader, we shall give a self-contained proof of this result even though the same argument has already appeared in an earlier paper of ours\footnote{G.~Trentinaglia, ``Tannaka duality for proper Lie groupoids,'' preprint. Posted on the ArXiv under a different title.} in a different theoretical framework.

\begin{teorema}\label{thm:SURJECT}
Let $\mathcal G$ be a proper Lie groupoid. Then the canonical homomorphism $\canhom[\mathcal G]: \mathcal G \longto \bidual{\mathcal G}$ is full, i.e., surjective.
\end{teorema}\begin{proof}
We start by proving that if the set $\mathcal G(x,x')$ is empty then the same is true of the set $\bidual{\mathcal G}(x,x')$.

Let $\varphi: {{\mathcal G x} \cup {\mathcal G x'}} \to \nC$ be the function which takes the value one on the orbit ${\mathcal G x}$ and the value zero on the orbit ${\mathcal G x'}$. (Recall that these are closed submanifolds.) This function is well-defined, because $\mathcal G(x,x')$ is empty. By Lemma \ref{lem:EQVEXT}, there is a global invariant smooth function $\Phi$ extending $\varphi$. Being invariant, $\Phi$ determines an endomorphism $a$ of the trivial representation $\underline\nC$ such that $a_z = {\Phi(z)\, \id_\nC}$ for all $z$; in particular, $a_x = \id$ and $a_{x'} = 0$. Now, suppose there is some $\lambda \in \bidual{\mathcal G}(x,x')$. Because of the naturality of $\lambda$, the existence of the morphism $a$ contradicts the invertibility of the linear map $\lambda(\underline\nC)$.

We are therefore reduced to proving that the induced isotropy group homomorphisms $(\canhom[\mathcal G])_x: \mathcal G_x \to \bidual{\mathcal G}_x$ are surjective for all $x$. So, let $x$ be an arbitrary base point. As $x$ will be fixed throughout the rest of the proof, put $\pi \equiv (\canhom[\mathcal G])_x$ and, for any representation $R = (E,\varrho)$ of $\mathcal G$, let $\pi_R$ denote the representation of the Lie group $\mathcal G_x$ on the vector space $E_x$ given by $g \mapsto \pi(g)(R)$. Also, let $\mathcal C$ denote the category \rep{\mathcal G}, and $F$ the functor \forget{\mathcal G,x} from $\mathcal C$ into the category of finite dimensional complex vector spaces (\ref{N.v2}).

Put $K \equiv \kernel\pi \subset \mathcal G_x$. This is a closed normal subgroup, because it coincides with the intersection ${\bigcap \kernel{\pi_R}}$ over all objects $R$ of $\mathcal C$. On the quotient $G \equiv \mathcal G_x/K$ there is a unique (compact) Lie group structure such that the quotient homomorphism $\mathcal G_x \to G$ becomes a Lie group homomorphism. Every $\pi_R$ can be indifferently thought of as a continuous representation of $\mathcal G_x$ or a continuous representation of $G$, and every linear map $A: F(R) \to F(S)$ is a morphism of $G$-modules if and only if it is a morphism of $\mathcal G_x$-modules. Being continuous, every $\pi_R$ is also smooth.

We claim there exists an object $R_0$ of $\mathcal C$ such that the corresponding $\pi_{R_0}$ is faithful as a representation of $G$. Indeed, by the compactness of the Lie group $G$, we can find $R_1, \ldots, R_\ell \in \Ob(\mathcal C)$ with the property that
\begin{equazione}\label{O.esp83}
{\kernel{\pi_{R_1}} \cap \cdots \cap \kernel{\pi_{R_\ell}}} = \{e\}\text,
\end{equazione}
where $e$ denotes the unit of $G$; compare \cite{BtD85}, p.~136. Then, if we set $R_0 = {R_1 \oplus \cdots \oplus R_\ell}$, the representation $\pi_{R_0}$ will be faithful because of the existence of an isomorphism of $G$-modules
\begin{equazione}\label{O.equ12}
F({R_1 \oplus \cdots \oplus R_\ell}) \iso {F(R_1) \oplus \cdots \oplus F(R_\ell)}\text.
\end{equazione}

Now, every irreducible, continuous, finite dimensional, complex $G$-module $V$ embeds as a submodule of some tensor power ${F(R_0)^{\otimes k} \otimes (F(R_0)^*)^{\otimes\ell}}$; see, for instance, \cite{BtD85}, p.~137. Since each $\pi(g)$ is a self-conjugate, tensor preserving natural transformation, this tensor power will be naturally isomorphic to $F\left({{R_0}^{\otimes k} \otimes ({R_0}^*)^{\otimes\ell}}\right)$ as a $G$-module and hence, for each object $V$ of the category \rep{G} of all continuous, finite dimensional, complex $G$-modules, there will be some object $R \in \Ob(\mathcal C)$ such that $V$ embeds into $F(R)$ as a submodule.

Next, consider an arbitrary natural transformation $\lambda \in \End{F}$. Let $R$ be an object of the category $\mathcal C$, and let $V \subset F(R)$ be a submodule. The choice of a complement to $V$ in $F(R)$ determines an endomorphism of modules $P_V: F(R) \to V \into F(R)$, which, by Proposition \ref{prop:EQVEXT}, must come from some endomorphism of $R$ in $\mathcal C$. This implies that the linear operators $\lambda(R)$ and $P_V$ on the space $F(R)$ commute with one another and, consequently, that $\lambda(R)$ maps the subspace $V$ into itself. Hereafter, we shall omit any reference to $R$ and write simply $\lambda_V$ for the linear map that $\lambda(R)$ induces on $V$ by restriction. Note finally that, given another submodule $W \subset F(S)$ and any equivariant map $B: V \to W$, one has the following identity of linear maps:
\begin{equazione}\label{O.equ11}
{B \circ \lambda_V} = {\lambda_W \circ B}\text.
\end{equazione}
To prove this identity, one first extends $B$ to an equivariant map $F(R) \to F(S)$ and then invokes Proposition \ref{prop:EQVEXT} once again.

Let $\mathsf F_G$ denote the functor $\rep{G} \longto \veb{}$ that assigns each $G$-module the underlying vector space. We will now define an isomorphism of complex algebras
\begin{equazione}\label{x.8}
\theta: \End{F} \isoto \End{\mathsf F_G}
\end{equazione}
so that the following diagram commutes
\begin{equazione}\label{O.equ13}\begin{split}
\xymatrix@C=35pt@R=20pt{\mathcal G_x \ar@{->>}[d] \ar[r]^-\pi & \End{F} \ar[d]_-\simeq^-\theta \\ G \ar[r]^-{\canhom[G]} & \End{\mathsf F_G}\text.\!\!}
\end{split}\end{equazione}
For each $G$-module $V$, there exists an object $R$ of $\mathcal C$ together with an embedding $V \into F(R)$, so we could define ${\theta(\lambda)}(V)$ as the restriction $\lambda_V$ of $\lambda(R)$ to $V$. Of course, it is necessary to check that this does not depend on the choices involved. Let two objects $R, S \in \Ob(\mathcal C)$ be given, along with two equivariant embeddings of $V$ into $F(R)$ and $F(S)$ respectively. Since it is always possible to embed everything equivariantly into $F({R\oplus S})$ without affecting the induced $\lambda_V$, it is no loss of generality to assume $R=S$. Now, from the identity (\ref{O.equ11}) above, it follows that the two embeddings actually determine the same linear endomorphism of $V$. This shows that $\theta$ is well-defined. The same identity also implies that $\theta(\lambda) \in \End{\mathsf F_G}$. On the other hand, put, for $\mu \in \End{\mathsf F_G}$ and $R \in \Ob(\mathcal C)$, $\mu^F(R) := \mu(F(R))$. Then $\mu^F \in \End{F}$ and $\theta(\mu^F) = \mu$, because of the existence of embeddings of the form $V \into F(R)$ and because of the naturality of $\mu$. This shows that $\theta$ is surjective, and also injective since $\lambda(R) = \theta(\lambda)(F(R))$. Finally, it is straightforward to check that (\ref{O.equ13}) commutes.

In order to conclude the proof, it will be enough to check that $\theta$ induces a bijection between $\mathrm{End}^{\overline\otimes}(F)$ and $\mathrm{End}^{\overline\otimes}(\mathsf F_G)$, for then our claim that $\pi$ is surjective will follow immediately from the commutativity of (\ref{O.equ13}) and the classical Tannaka duality theorem for compact groups (which says that \canhom[G] establishes a bijection between $G$ and $\mathrm{End}^{\overline\otimes}(\mathsf F_G)$; see, for example, \cite{JoSt91} or \cite{BtD85}). This can safely be left to the reader.
\end{proof}

\begin{definizione}\label{def:REFLEX}
We call a Lie groupoid $\mathcal G$ \textit{reflexive,} when the corresponding canonical homomorphism \canhom[\mathcal G] is an isomorphism of $\C^\infty$-groupoids.
\end{definizione}

\begin{teorema}\label{thm:INJ=>HOMEO}
Suppose $\mathcal G$ is proper. Then, the canonical homomorphism $\canhom[\mathcal G]: \mathcal G \longto \bidual{\mathcal G}$ is an isomorphism of topological groupoids if, and only if, it is faithful (i.e., injective).
\end{teorema}\begin{proof}
The only thing we have not yet proved is that when \canhom[\mathcal G] is faithful, for any open subset $\Gamma$ (of the manifold of arrows) of the groupoid $\mathcal G$ and for any point $g \in \Gamma$, the image $\canhom[\mathcal G](\Gamma)$ is a neighbourhood of $\canhom[\mathcal G](g)$ in (the space of arrows of) the groupoid \bidual{\mathcal G}.

We start by showing that there exists a representation $R = (E,\varrho)$ of~$\, \mathcal G$, and an open ball $P \subset \GL{E}$ centred at $\varrho(g)$, such that $\varrho^{-1}(P) \subset \Gamma$. Suppose $g \in \mathcal G(x,x')$. By one of the remarks we made in the course of the proof of Theorem \ref{thm:SURJECT}, there is a representation $R = (E,\varrho)$ whose associated isotropy homomorphism $\varrho_x: \mathcal G_x \to \Aut{E_x}$ is faithful. The restriction of such a representation to the subset $\mathcal G(x,x')$ is injective. Fix a decreasing sequence
$$%
\cdots \subset P_{i+1} \subset P_i \subset \cdots \subset P_1
$$%
of open balls in \GL{E} converging to $\varrho(g)$. Put
$$%
\Sigma_i := {\varrho^{-1}\bigl(\overline{P_i}\bigr) - \Gamma}\text.
$$%
The sets $\Sigma_i$ form a decreasing sequence of closed subsets of the manifold of arrows of $\mathcal G$, with empty intersection. By the properness of $\mathcal G$, there exists some $i$ such that $\Sigma_i = \varnothing$. This proves the claim.

With this information at hand, we can conclude at once by the surjectivity of \canhom[\mathcal G] (Theorem \ref{thm:SURJECT}). Indeed, on the one hand, we have
$$%
\canhom[\mathcal G](g)\, \in\, {\eval[R]}^{-1}(P)\, =\: {\canhom[\mathcal G]\: {\canhom[\mathcal G]}^{-1}\, {\eval[R]}^{-1}(P)}\, =\, \canhom[\mathcal G]\bigl(\varrho^{-1}(P)\bigr)\, \subset\, \canhom[\mathcal G](\Gamma)\text.
$$%
On the other hand, ${\eval[R]}^{-1}(P)$ is already known to be open.
\end{proof}

The following lemma is a direct consequence of a more general statement which will be proved in the next section.

\begin{lemma}\label{lem:INJREP}
Let $\mathcal G$ be a Lie groupoid, and let $\varrho: \mathcal G \to \GL{E}$ be a representation of~$\,\mathcal G$. Let $\Gamma$ be an open subset (of the manifold of arrows) of~$\,\mathcal G$. Suppose that for some pair of base points $x,x'$, the representation $\varrho$ restricts to an injection on the subset $\Gamma(x,x') := {\Gamma \cap \mathcal G(x,x')}$. Then $\varrho$ is an immersion on some open neighbourhood of the same subset.
\end{lemma}\begin{proof}
This is the special case of Lemma \ref{lem:INJ=>IMM} (see the next section) where $\mathcal H = \mathcal G$, $\mathcal H' = \mathcal G' = \Gamma' = \GL{E}$, $\rho = \rho' =$ identity, and $\psi = \varrho$.
\end{proof}

We shall say that a Lie groupoid $\mathcal G$ \textit{has enough representations,} when for each base point $x$ and for each $g \in \mathcal G(x,x')$ with $g \neq x$, there exists some representation $\varrho: \mathcal G \to \GL{E}$ such that $\varrho(g) \neq \id_{E_x}$. Then, we have the following characterization of reflexivity for proper Lie groupoids:

\begin{teorema}\label{thm:REFLEX}
A proper Lie groupoid is reflexive if and only if it possesses enough representations.
\end{teorema}\begin{proof}
The existence of enough representations implies that the canonical homomorphism is faithful and hence that it induces, in view of Theorem \ref{thm:INJ=>HOMEO}, a {\em homeomorphism} between the (space of arrows of the) given proper Lie groupoid $\mathcal G$ and (that of) the corresponding Tannakian bidual \bidual{\mathcal G}. Thus, it will suffice to show that \canhom[\mathcal G] is a {\em local isomorphism} of $\C^\infty$-spaces.

Let $g \in \mathcal G(x,x')$ be given. As it was already observed in the proof of Theorem \ref{thm:INJ=>HOMEO}, we can find $R = (E,\varrho) \in \Ob\:\rep{\mathcal G}$ such that the restriction $\varrho_{x,x'}: \mathcal G(x,x') \to \Iso{E_x}{E_{x'}}$ is injective. Then, by Lemma \ref{lem:INJREP}, $\varrho$ induces a diffeomorphism $\varrho|_\Gamma: \Gamma \simeq \varrho(\Gamma)$ between an open neighbourhood $\Gamma$ of $g$ in $\mathcal G$ and a submanifold $\varrho(\Gamma)$ of the linear groupoid \GL{E}. Now, put $\Omega \equiv \canhom[\mathcal G](\Gamma)$. We know that $\canhom[\mathcal G]|_\Gamma: \Gamma \to \Omega$ is a homeomorphism and a $\C^\infty$-mapping. Moreover, we know that the evaluation representation $\eval[R]|_\Omega: \Omega \to \GL{E}$ induces a $\C^\infty$-mapping from $\Omega$ onto $\eval_R(\Omega) = \varrho(\Gamma)$. Hence, $(\canhom[\mathcal G]|_\Gamma)^{-1} = {(\varrho|_\Gamma)^{-1} \circ \eval[R]|_\Omega}$ is also a $\C^\infty$-mapping, from $\Omega$ onto $\Gamma$. This finishes the proof.
\end{proof}

This is a good point to discuss some examples. Actually, we start with a counterexample \cite{LuOl01,2008}, which shows that not every proper Lie groupoid is reflexive. This counterexample makes it clear that the classical duality theorem of Tannaka cannot be generalized to proper Lie groupoids, at least in the framework of representations on vector bundles.\footnote{One can define a notion of representation on a more general type of linear bundle which makes this generalization possible, though. Compare \cite{2008}.} Recall that a \textit{bundle of} (\textit{compact}) \textit{Lie groups} is a (proper) Lie groupoid whose source and target map coincide.

\begin{paragrafetto}[Example: a bundle of compact Lie groups not having enough representations]\label{xmp:CTRXMP}
We start by constructing, for any Lie group $H$ and any automorphism $\chi \in \Aut{H}$, a (locally trivial) bundle of Lie groups with fibre $H$ over the unit circle $\mathrm S^1$, hereafter denoted by $H[\chi]$.

As a topological space, $H[\chi]$ is defined to be the quotient of the Cartesian product ${\nR \times H}$ induced by the equivalence relation
\begin{equazione}\label{def:H[chi]"1}
(t,h) \thicksim (t',h') \quad \aeq \quad t'-t \in \nZ \quad \text{and} \quad h' = \chi^{t'-t}(h)\text.
\end{equazione}
The bundle fibration $H[\chi] \to \mathrm S^1$ is induced by the map $(t,h) \mapsto \exp(2\pi i t)$. In terms of representatives of equivalence classes, the composition law is given by
\begin{equazione+}\label{def:H[chi]"2}
{[t',h'] \cdot [t,h]} := \bigl[t',{h' \chi^{t'-t}(h)}\bigr]\text. & \bigl(\text{Defined when $t'-t \in \nZ$.}\bigr)
\end{equazione+}
As to the differentiable structure, consider the open cover of $\mathrm S^1$ given by the two open intervals $I = (0,1)$ and $J = \bigl(-\tfrac12,\tfrac12\bigr)$ (more precisely, the images thereof under the exponential map). One has two corresponding smoothly compatible local trivializing charts for $H[\chi] \to \mathrm S^1$, namely
\begin{equazione}\label{def:H[chi]"3}
\left\{\begin{aligned}
H[\chi]_I \isoto {I \times H}\text, \quad &\text{given by} \quad [t,h] \mapsto (\exp(2\pi i t),h)\text{, $0<t<1$,} \quad \text{and}
\\
H[\chi]_J \isoto {J \times H}\text, \quad &\text{given by} \quad [t,h] \mapsto (\exp(2\pi i t),h)\text{, $-\tfrac12<t<\tfrac12$.}
\end{aligned}\right.
\end{equazione}
These charts determine the differentiable structure. The transition mapping is given by the identity over $(0,1/2)$ and by the diffeomorphism $(z,h) \mapsto \bigl(z,\chi(h)\bigr)$ over $(1/2,1)$. So $H[\chi]$ will not be globally trivial, in general.

We proceed to consider the following special case. Take $H := \mathrm T^2$ to be the two-torus. Fix an arbitrary integer $\ell \in \nZ$, and consider the map
\begin{equazione}\label{def:chi-ell}
\chi_\ell: \mathrm T^2 \to \mathrm T^2 \quad \text{given by} \quad (a,b) \mapsto (a,a^\ell b)\text;
\end{equazione}
$\chi_\ell$ is an automorphism of the Lie group $\mathrm T^2$ (its inverse being $\chi_{-\ell}$).

It is a nice ex\-er\-cise\inciso{involving only some elementary knowledge of the representation theory of Lie groups}to show that (provided $\ell \neq 0$) the image of the following embedding of bundles of Lie groups over $\mathrm S^1$
\begin{equazione}\label{def:phi-ell}
\phi_\ell: {\mathrm S^1 \times \mathrm T^1} \into \mathrm T^2[\chi_\ell]\text, \quad (\exp(2\pi i t),b) \mapsto [t,(1,b)]
\end{equazione}
is contained in the kernel of every representation of the Lie groupoid $H[\chi]$. In fact, one can say more. Consider the following map
\begin{equazione}\label{def:psi-ell}
\psi_\ell: \mathrm T^2[\chi_\ell] \onto {\mathrm S^1 \times \mathrm T^1}\text, \quad [t,(a,b)] \mapsto (\exp(2\pi i t),a)\text.
\end{equazione}
This is an epimorphism of bundles of Lie groups over $\mathrm S^1$. Its kernel is precisely the image of the embedding $\phi_\ell$. Then, it is easy to show that the pullback along $\psi_\ell$ yields an isomorphism of representation categories
\begin{equazione}\label{O.equ10}
{\psi_\ell}^*: \underline{\smash[b]{\mathsf{Rep}}}\bigl(\mathrm S^1 \times \mathrm T^1\bigr) \xto{\quad \simeq \quad} \underline{\smash[b]{\mathsf{Rep}}}\bigl(\mathrm T^2[\chi_\ell]\bigr)
\end{equazione}
{\em over the base manifold~$\,\mathrm S^1$,} i.e., compatible with the two canonical forgetful functors into \veb{\mathrm S^1}. Since two isomorphic (or even equivalent) representation categories over the same base manifold must be regarded, for all purposes of any reasonable duality theory, as indistinguishable, one is led to the conclusion that there cannot be any procedure by means of which one can possibly recover the groupoids $\mathrm T^2[\chi_\ell]$, $\ell \neq 0$.
\end{paragrafetto}

After this counterexample to reflexivity, let us turn our attention to some positive examples.

\begin{paragrafetto}[Example: transitive proper Lie groupoids]\label{xmp:TRANSLGPD}
Recall that a Lie groupoid $\mathcal G$ is said to be \textit{transitive} when the combined source--target map
\begin{equazione}\label{def:TRANSLGPD}
(\s,\t): \mca[1]{\mathcal G} \longto {\mca[0]{\mathcal G} \times \mca[0]{\mathcal G}}\text, \qquad g \mapsto (\s(g),\t(g))
\end{equazione}
is a {\em surjective submersion.} Any transitive proper Lie groupoid is Morita equivalent to a compact Lie group; see, for instance, \cite{MoeMrc03}. It is an easy exercise, in view of the remarks about Morita equivalences at the end of the preceding section, to show that any such groupoid is reflexive.
\end{paragrafetto}

\begin{paragrafetto}[Example: action groupoids associated with compact Lie group actions]\label{xmp:CPTLGPACT}
For each smooth (left) action of a Lie group $G$ on a manifold $M$, there is a Lie groupoid ${G \ltimes M}$ with base $M$, called the \textit{action} (or \textit{translation}) \textit{groupoid associated with the given action of~$\,G$ on $M$}. One takes the Cartesian product ${G \times M}$ as the manifold of arrows, the projection $(g,x) \mapsto x$ as the source, the group action $(g,x) \mapsto gx$ as the target, and the operation
\begin{equazione}\label{def:action-gpd}
(g',x')(g,x) = (g'g,x)
\end{equazione}
as the composition law. Any action groupoid ${K \ltimes M}$ associated with a smooth action of a {\em compact} Lie group $K$ on a manifold $M$ is reflexive. Indeed, if $\boldsymbol V$ is any faithful $K$-module, in other words, any faithful representation $K \into \GL{\boldsymbol V}$ on a finite dimensional vector space $\boldsymbol V$, then we get a corresponding faithful representation of the groupoid ${K \ltimes M}$ on the trivial vector bundle ${M \times \boldsymbol V}$:
\begin{equazione}\label{def:action-rep}
(k,x)\: \mapsto\: \bigl\{ (x,\boldsymbol v) \mapsto (kx,{k\cdot \boldsymbol v}) \bigr\}\text.
\end{equazione}
\end{paragrafetto}


\begin{paragrafetto}[Example: source-proper \'etale Lie groupoids]\label{xmp:source-proper+'etale}
Recall that a Lie groupoid is said to be \textit{\'etale,} if its source and target map are \'etale, that is, local smooth isomorphisms. We shall say that a Lie groupoid is \textit{source-proper,} whenever the corresponding source map is proper.

For any source-proper \'etale Lie groupoid $\mathcal G$, one has the corresponding \textit{regular representation $(R,\varrho) = (R^{\mathcal G},\varrho^{\mathcal G}) \in \Ob\:\rep{\mathcal G}$}, defined as follows. For each base point $x$, let $\ell(x) := \norm{\s^{-1}(x)}$ denote the cardinality of the respective source fibre (which is a finite set, as a consequence of the fact that, by our hypotheses on $\mathcal G$, the source must be a locally trivial map with discrete fibres). Then put
\begin{equazione}\label{def:R_x}
R_x := \C^0\bigl(\s^{-1}(x);\nC\bigr) \iso \nC^{\ell(x)}\text.
\end{equazione}
The local identifications $R_U \iso {U \times \nC^\ell}$ obtained in the obvious way from the local trivializations $\mathcal G^U \iso {U \times \{1, \ldots, \ell\}}$ for the source map of $\mathcal G$ provide a convenient atlas of local trivializing charts for the vector bundle $R$. Define
\begin{equazione+}\label{def:regrep}
\varrho(g): R_x \to R_{x'} \quad \text{as} \quad f \mapsto {\varrho(g)}(f) := f(\text-\,g)\text. & \bigl(\text{where $g \in \mathcal G(x,x')$.}\bigr)
\end{equazione+}
This action of~$\,\mathcal G$ on $R$ is smooth\inciso{because in any local trivializing charts it looks like a constant permutation}and, of course, {\em faithful.}
\end{paragrafetto}

\begin{paragrafetto}[Example: effective orbifolds]\label{xmp:EFFORB}
For any smooth manifold $M$, let ${\Gamma M}$ denote the groupoid over $M$ whose arrows $x \to x'$ are the germs of local diffeomorphisms in $M$ mapping $x$ to $x'$. If $\mathcal G$ is an \'etale Lie groupoid over $M$, one has a canonical homomorphism of groupoids $\mathcal G \to {\Gamma M}$ over $M$, the so-called \textit{effect of~$\,\mathcal G$}, defined by sending any arrow $g$ to the germ of local smooth isomorphism associated with a sufficiently small open neighborhood of $g$ in the manifold of arrows of $\mathcal G$. An \textit{effective} Lie groupoid is, by definition, an \'etale Lie groupoid whose effect is faithful. Compare \cite{MoeMrc03}.

For any \'etale Lie groupoid, one has a representation on the tangent bundle of the respective base manifold, the \textit{tangent representation,} obtained by assigning every arrow $g$ the differential of the effect of $g$ at the point $\s(g)$. In the proper, effective case, the tangent representation can easily be shown to be {\em faithful.} For a proof of this fact, we refer the reader to Section 28 of \cite{2008}.
\end{paragrafetto}

\sezione{Representative charts}\label{sez:REPCHARTS}
The counterexample to reflexivity presented in the previous section (Example \ref{xmp:CTRXMP}) motivates the following

\begin{definizione}\label{def:parareflexive}
We shall say that a Lie groupoid $\mathcal G$ is \textit{parareflexive,} if the Tannakian bidual \bidual{\mathcal G} is a {\em Lie} groupoid and the canonical homomorphism $\canhom[\mathcal G]: \mathcal G \to \bidual{\mathcal G}$ is a {\em surjective submersion.}
\end{definizione}

Let us point out a couple of direct consequences of this definition. First of all, the categories of representations of any parareflexive Lie groupoid and of the corresponding Tannakian bidual are the same, in view of the following result (which we anticipated in Section \ref{sez:BIDUAL}):

\begin{proposizione}\label{prop:EVALFUNC}
Let $\mathcal G$ be a Lie groupoid whose associated canonical homomorphism \canhom[\mathcal G] is full, i.e., surjective. Then the evaluation functor
$$%
\eval: \rep{\mathcal G} \longto \rep{\bidual{\mathcal G}}\text, \qquad R = (E,\varrho) \mapsto (E,\eval[R])
$$%
is an isomorphism, having the pullback ${\canhom[\mathcal G]}^*$ for inverse.
\end{proposizione}\begin{proof}Trivial.\end{proof}

Secondly, the bidual of any parareflexive Lie groupoid must be a reflexive Lie groupoid. This follows from Proposition \ref{prop:EVALFUNC} and the remark that if
$$%
\xymatrix@C=40pt@R=20pt{\rep{\mathcal G} \ar[d]^{\forget{\mathcal G}} \ar[r]^-\simeq_-{\boldsymbol\Phi} & \rep{\mathcal H} \ar[d]^{\forget{\mathcal H}} \\ \veb{M} \ar@{=}[r] & \veb{M}}
$$%
is an isomorphism (more generally, an equivalence) between the representation categories of two Lie groupoids over the same base manifold $M$, then there is a corresponding $\C^\infty$-isomorphism $\phi: \bidual{\mathcal H} \to \bidual{\mathcal G}$, which, in the special case $\mathcal H \equiv \bidual{\mathcal G}$, $\boldsymbol\Phi \equiv \eval$, will provide a $\C^\infty$-inverse for \canhom[\bidual{\mathcal G}]. (Compare \cite{2008}, \S24.)

\medskip

In the present and in the next section, we shall delve into the property of parareflexivity for proper Lie groupoids. It is our purpose to obtain, for such groupoids, an explicit necessary and sufficient condition for parareflexivity, hopefully a direct generalization of the characterization of reflexivity we obtained in the preceding section (Theorem \ref{thm:REFLEX}).

To motivate the next definition, let us consider any parareflexive, proper Lie groupoid $\mathcal G$. By the properness of $\mathcal G$, we know (from the proof of Thm.\ \ref{thm:SURJECT}, for instance) that for any pair $x,x'$ of base points there exists a representation $R = (E,\varrho)$ of $\mathcal G$ such that the corresponding evaluation representation $\eval[R]: \bidual{\mathcal G} \to \GL{E}$ restricts to an injection on the subset $\bidual{\mathcal G}(x,x')$. Under the assumption that $\mathcal G$ is parareflexive, Lemma \ref{lem:INJREP} implies that the smooth representation \eval[R] is an immersion in the vicinity of the same subset. Thus, any given arrow of the bidual \bidual{\mathcal G} will have an open neighbourhood $\Omega$ (in the manifold of arrows of the same groupoid) such that \eval[R] induces a dif\-feo\-morph\-ism\inciso{a fortiori, a bi\-jec\-tion}onto a submanifold $R(\Omega) := \eval[R](\Omega)$ of the manifold of arrows of the linear groupoid \GL{E}.

\begin{definizione}\label{def:REPCHART}
Let $\mathcal G$ be a proper Lie groupoid. We shall call \textit{representative chart} (\textit{for the bidual of~$\,\mathcal G$}) any pair $(\Omega,R)$ consisting of an open subset $\Omega$ of the space of arrows of the groupoid \bidual{\mathcal G} and an object $R = (E,\varrho) \in \Ob\:\rep{\mathcal G}$ such that $\eval_R: \bidual{\mathcal G} \to \GL{E}$ induces a bijective correspondence between $\Omega$ and a submanifold of the manifold of arrows of the linear groupoid $\GL{E}$.
\end{definizione}

\begin{paragrafetto}[Remark]\label{rmk:def:REPCHART}
This notion is tailor-made for {\em proper} Lie groupoids. By suitably modifying it, one can generalize some of the results below\inciso{Prop.\ \ref{prop:CRIT"1}, for instance}to arbitrary Lie groupoids. We shall refrain from doing this here. For the general theory, we refer the reader to Chapter V of \cite{2008}.
\end{paragrafetto}

In order to make Definition \ref{def:REPCHART} more plausible, we proceed to show that, in fact, for any representative chart $(\Omega,R)$ the bijection $\eval[R]: \Omega \isoto R(\Omega) := \eval[R](\Omega)$ is a {\em homeomorphism.} From the considerations preceding Definition \ref{def:REPCHART}, it will follow that when $\mathcal G$ is parareflexive, the map $\lambda \mapsto \lambda(R)$ must induce a {\em diffeomorphism} between $\Omega$ and a submanifold $R(\Omega)$ of \GL{E}. This justifies the name of ``chart''.

\begin{lemma}\label{lem:Gamma->Sigma}
Suppose that $\mathcal G$ is proper. Let $(E,\varrho)$ be a representation of $\mathcal G$, and let $\Gamma$ be an open subset of (the manifold of arrows of) $\mathcal G$ such that the image $\Sigma \equiv \varrho(\Gamma)$ is a submanifold of (the manifold of arrows of) \GL{E}. Then the restriction of $\varrho$ to $\Gamma$ is a submersion onto $\Sigma$.
\end{lemma}\begin{proof}
To begin with, we observe that the map $\Sigma \to M$ induced by the source map of \GL{E} by restriction is a submersion. Indeed, for each $\sigma \in \Sigma$, let us say $\sigma = \varrho(\gamma)$ with $\gamma \in \Gamma$, there exists some smooth local source section $U \to \Gamma$, defined over a neighbourhood $U$ of the base point $x \equiv \s(\gamma)$, sending $x \mapsto \gamma$; hence, upon composing with $\varrho$, we obtain also a smooth local source section $U \to \Sigma$ sending $x \mapsto \sigma$.

Fix now $\gamma \in \Gamma$, and put $\sigma \equiv \varrho(\gamma)$. Let $x \equiv \s(\gamma)$. By the preceding remark, there exist, in some open neighbourhoods $\Gamma(\gamma) \subset \Gamma$ of $\gamma$ and $\Sigma(\sigma) \subset \Sigma$ of $\sigma$, trivializing charts for the source map (on $\Gamma$ and $\Sigma$ respectively). In suitably chosen such charts, the restriction of $\varrho$ will take the following form:
$$%
{U \times \nR^m} \iso \Gamma(\gamma) \xto{\quad \varrho|_{\Gamma(\gamma)} \quad} \Sigma(\sigma) \iso {U \times \nR^n}\text, \qquad (u,\boldsymbol x) \mapsto \bigl(u,\boldsymbol y(u,\boldsymbol x)\bigr)
$$%
(for some open neighbourhood $U$ of $x$ in $M$). If we can show that the partial map $\boldsymbol x \mapsto \boldsymbol y(x,\boldsymbol x)$ is submersive at zero, we are clearly done. To this end, we proceed to check that the restriction
$$%
\varrho|_{\mathcal G(x,\text-)}: \mathcal G(x,\text-) \xto\qquad \varrho\bigl(\mathcal G(x,\text-)\bigr)
$$%
is a submersion (at $\gamma$ now generic) onto a submanifold of \GL{E}.

Choose any open subset $U' \subset M$, with $x' \equiv \t(\gamma) \in U'$, small enough to ensure that a local $G$-equivariant chart $\mathcal G(x,S) \iso {S \times G}$ (where $G$ denotes the isotropy group of $\mathcal G$ at $x$) can be found over $S \equiv {(\mathcal G\cdot x) \cap U'}$. It is no loss of generality to assume that $\gamma$ corresponds to $(x',e)$ in this chart ($e$ being the unit of $G$). We obtain a smooth section $s \mapsto (s,e) \stackrel{\iso}{\mapsto} g \mapsto \varrho(g)$ to the target map of $\GL{E}$ over $S$. Next, observe that the isotropy homomorphism $\varrho_x: G \to \GL{E_x}$ induced by $\varrho$ factors through the quotient $K \equiv G/\kernel{\varrho_x}$. The closed Lie subgroup $K \into \GL{E_x}$ and the target section $S \into \GL{E}$ can be combined into an embedding ${S \times K} \into \GL E$ closing the diagram
$$%
\xymatrix@C=30pt@R=20pt{{S \times G}\ar@{->>}[d]\ar[r]^-\iso & \mathcal G(x,S)\ar[d]^\varrho \\ {\;S \times K\;}\ar@{^(-->}[r] & \GL{E}\text.\!\!}
$$%
The proof is now evidently finished, because $\mathcal G(x,S) = {\mathcal G(x,\text-) \cap \t^{-1}(U')}$ and $\varrho\bigl(\mathcal G(x,S)\bigr) = {\varrho\bigl(\mathcal G(x,\text-)\bigr) \cap \t^{-1}(U')}$.
\end{proof}

Now, if we make $\Gamma \equiv (\canhom[\mathcal G])^{-1}(\Omega)$, we have $\varrho(\Gamma) = R(\Omega)$, because of the identity ${\eval[R] \circ \canhom[\mathcal G]} = \varrho$ and because of the surjectivity of \canhom[\mathcal G] (Theorem \ref{thm:SURJECT}). On the other hand, Lemma \ref{lem:Gamma->Sigma} tells us that $\varrho: \Gamma \to R(\Omega)$ is a submersion. The bijection $\eval_R: \Omega \isoto R(\Omega)$ is therefore an open mapping and hence, being also continuous, it must be a homeomorphism.

\begin{paragrafetto}[Remarks]\label{rmk:REPCHARTS}
{\bfseries (a)~}~Given two isomorphic objects $R \iso S$ in the category \rep{\mathcal G}, the pair $(\Omega,R)$ is a representative chart for the bidual of $\mathcal G$ if, and only if, the same is true of the pair $(\Omega,S)$.

{\bfseries (b)~}~Let $(\Omega,R)$ be a representative chart for the bidual of $\mathcal G$. Then, for each open subset $\Omega' \subset \Omega$, the same is true of the pair $(\Omega',R)$. (This is a consequence of the result that $\eval_R: \Omega \to R(\Omega)$ must be a homeomorphism.)
\end{paragrafetto}

A first step in the direction of an effective characterization of parareflexivity in the proper case is represented by the following

\begin{proposizione}\label{prop:CRIT"1}
Let $\mathcal G$ be a proper Lie groupoid. Then the Tannakian bidual \bidual{\mathcal G} is a Lie groupoid if, and only if, the following two conditions are satisfied:
\begin{elenco}
\item the domains of representative charts cover the space of arrows of the bidual groupoid \bidual{\mathcal G}, i.e., for each arrow $\lambda$ in \bidual{\mathcal G} there exists a representative chart $(\Omega,R)$ with $\lambda \in \Omega$;
\item if~$\,(\Omega,R)$ is a representative chart, the same is true of~$\,(\Omega,{R \oplus S})$ for any other representation $S$ of~$\,\mathcal G$.
\end{elenco}
\end{proposizione}

Before we embark into the proof of this proposition, we need to establish a technical lemma (compare Lemma \ref{lem:INJREP}).

\begin{lemma}\label{lem:INJ=>IMM}
Let $\rho: \mathcal G \to \mathcal H$, $\rho': \mathcal G' \to \mathcal H'$, and $\psi: \mathcal H \to \mathcal H'$ be homomorphisms of Lie groupoids over a manifold $M$. (It is understood that they are identical on $M$.) Let $\Gamma$, resp., $\Gamma'$ be an open subset of (the manifold of arrows of) $\mathcal G$, resp., $\mathcal G'$ such that the image $\Sigma \equiv \rho(\Gamma)$, resp., $\Sigma' \equiv \rho'(\Gamma')$ is a submanifold of (the manifold of arrows of) $\mathcal H$, resp., $\mathcal H'$. Suppose, in addition, that the induced maps $\Gamma \to \Sigma$ and~$\,\Gamma' \to \Sigma'$ are open.

Assume that $\psi$ sends $\Sigma$ into $\Sigma'$, and let $\psi_\Sigma: \Sigma \into \Sigma'$ be the induced map. Put $\Sigma(x,y) := {\Sigma \cap \mathcal H(x,y)}$, and let $\sigma \in \Sigma(x,y)$. Suppose finally that $\psi$ is injective on $\Sigma(x,y)$. Then $\psi_\Sigma$ is an immersion at $\sigma$.
\end{lemma}\begin{proof}
By reasoning as we did in the first half of the proof of Lemma \ref{lem:Gamma->Sigma}, we are reduced to showing that the restriction
\begin{equazione+}
\psi|_{\Sigma(x,\text-)}: \Sigma(x,\text-) \longto \Sigma'(x,\text-) & (\text{where~}\: \Sigma(x,\text-) := {\Sigma \cap \mathcal H(x,\text-)}\: \text{~etc.})
\end{equazione+}
is an immersion at $\sigma$ (note that $\Sigma(x,\text-) \subset \Sigma$ and $\Sigma'(x,\text-) \subset \Sigma'$ are submanifolds, because, as we observed at the beginning of the proof of Lemma \ref{lem:Gamma->Sigma}, the restriction of the source map of~$\,\mathcal H$, resp., $\mathcal H'$ to $\Sigma$, resp., $\Sigma'$ is a submersion).

The Lie groupoid homomorphism $\rho$ induces a homomorphism $\rho_x: \mathcal G_x \to \mathcal H_x$ between the isotropy groups at $x$. By factoring out the kernel of $\rho_x$, we obtain a {\em Lie subgroup $\eta: G \into \mathcal H_x$}, where $G \equiv \mathcal G_x/\kernel{\rho_x}$. It is not hard to see that, from our assumptions, it follows that there exist an open neighbourhood $A$ of the unit $e$ in $G$ and a smooth target section $\tau: Z \to \mathcal H(x,\text-)$ defined over a submanifold $Z \subset M$ containing $y$ with $\tau(y) = \sigma$, such that the mapping
\begin{equazione}\label{def:tau*eta}
{\tau \cdot \eta}: {Z \times A} \into \mathcal H(x,\text-)\text, \qquad (z,a) \mapsto {\tau(z) \cdot \eta(a)}
\end{equazione}
defines an open embedding of ${Z \times A}$ onto some open neighbourhood of $\sigma$ in $\Sigma(x,\text-)$, in other words, a local parametrization for $\Sigma(x,\text-)$ in the vicinity of $\sigma$ by the ``parameter set'' ${Z \times A}$. (More details can be found in \cite{2008}, \S22.)

With respect to any two (conveniently chosen) local parametrizations of the form (\ref{def:tau*eta}), let us say, ${\tau \cdot \eta}: {Z \times A} \into \Sigma(x,\text-)$ in the vicinity of $\sigma$ and ${\tau' \cdot \eta'}: {Z' \times A'} \into \Sigma'(x,\text-)$ in the vicinity of $\sigma'$, the restriction of $\psi$ reads
\begin{equazione}
{Z \times A} \to {Z' \times A'}\text, \qquad (z,a) \mapsto (z,a'(z,a))\text.
\end{equazione}
Thus, we are further reduced to showing that the partial map $a \mapsto a'(y,a)$ on $A$ is an immersion at the unit $e \in G$. (Note that this map is injective, by assumption.) Consider the following commutative diagram:
$$%
\xymatrix@C=38pt{\mathcal H_x \ar@{=}[r] & \mathcal H_x \ar[r]^-{\psi_x} & {\mathcal H'}_x \ar@{=}[r] & {\mathcal H'}_x \\ & \mathcal H(x,y) \ar[u]_{\tau(y)^{-1} \cdot} \ar[r]^-{\psi_{x,y}} & \mathcal H'(x,y) \ar[u]_{\psi(\tau(y))^{-1} \cdot} & \\ A \ar@{^(->}[uu]^\eta \ar@{=}[r] & {\{y\} \times A} \ar@{_(->}[u]_{\tau \cdot \eta} \ar[r] & {\{y\} \times A'} \ar@{_(->}[u]_{\tau' \cdot \eta'} \ar@{=}[r] & A'\ar@{^(->}[uu]^{\eta'}}
$$%
The bottom map in this diagram is precisely $a \mapsto a'(y,a)$. Now, it is not difficult to see that there is a unique Lie group homomorphism $\phi: G \to G'$ such that ${\eta' \circ \phi} = {\psi_x \circ \eta}$. Since $\phi$ agrees with $a \mapsto a'(y,a)$ on $A$, it must be immersive, because so must be any Lie group homomorphism which is injective in a neighbourhood of $e$. Hence $a \mapsto a'(y,a)$ is immersive at $e$.
\end{proof}

\begin{corollario}\label{cor:HOMEO=>DIFFEO}
Under the same assumptions as in Lemma \ref{lem:INJ=>IMM}, suppose, in addition, that $\psi_\Sigma$ is a homeomorphism of~$\,\Sigma$ onto $\Sigma'$. Then the same map is actually a diffeomorphism between $\Sigma$ and $\Sigma'$.\qed
\end{corollario}

\begin{proof}[Proof of Proposition \ref{prop:CRIT"1}]
(Necessity.)~~Let \bidual{\mathcal G} be a Lie groupoid. The assertion that Condition (i) is satisfied is actually the content of the remarks preceding Definition \ref{def:REPCHART}.

As to Condition (ii), let any representative chart $(\Omega,R)$ be given, with $R = (E,\varrho)$, and let $S = (F,\varsigma)$ be an arbitrary representation of $\mathcal G$. To begin with, we observe that the evaluation representation \eval[R\oplus S] injects $\Omega$ into \GL{E\oplus F}; indeed, since $\lambda(R\oplus S) = {\lambda(R) \oplus \lambda(S)}$ for all $\lambda$ in \bidual{\mathcal G}, we have that \eval[R\oplus S] factors through the submanifold (in fact, embedded subgroupoid)
\begin{equazione}\label{equ:GL(E)xGL(F)->GL(E+F)}
{\GL{E} \times_{M\times M} \GL{F}} \into \GL{E\oplus F}\text, \qquad (\mu,\nu) \mapsto {\mu \oplus \nu}
\end{equazione}
via the map $\lambda \mapsto \bigl(\eval[R](\lambda),\eval[S](\lambda)\bigr)$, which must be injective on $\Omega$, because, by hypothesis, so is \eval[R]. From this observation and from Lemma \ref{lem:INJ=>IMM}, it follows that \eval[R\oplus S] is an {\em immersion} on $\Omega$. We contend that \eval[R\oplus S] actually induces a {\em homeomorphism} between $\Omega$ and a subspace of \GL{E\oplus F}; once this is proved, it will follow that $({R\oplus S})(\Omega)$ is a submanifold of \GL{E\oplus F}. Now, let an open subset $\Omega' \subset \Omega$ be given. Fix an arbitrary open subset $\Lambda' \subset \GL{E}$ with ${R(\Omega) \cap \Lambda'} = R(\Omega')$ (such subsets exist because $\Omega$ and $R(\Omega)$ are homeomorphic via $\eval_R$). Then our contention follows at once from the identity
$$%
{({R\oplus S})(\Omega) \cap \bigl({\Lambda' \times_{M\times M} \GL F}\bigr)} = ({R\oplus S})(\Omega')\text.
$$%

(Sufficiency.)~~Suppose the two conditions (i) and (ii) are satisfied. We claim that when $(\Omega,R)$ is an arbitrary representative chart, the map \eval[R] establishes a $\C^\infty$-isomorphism between the $\C^\infty$-space $(\Omega,{\R^\infty|_\Omega})$ (see Section \ref{sez:BIDUAL}) and the submanifold $X \equiv R(\Omega)$ of the linear groupoid \GL{E}.

We know that $\eval[R]|_\Omega: \Omega \to X$ is a $\C^\infty$-mapping, because the evaluation representation $\eval[R]: \bidual{\mathcal G} \to \GL{E}$ is a $\C^\infty$-homomorphism. (Lemma \ref{lem:UNIVPROP}.)

The smoothness of the inverse bijection $(\eval[R]|_\Omega)^{-1}: X \to \Omega$ is less obvious. We have already proved that this is a homeomorphism. Now, recall our notation (\ref{N.v8}, \ref{N.v9}) from Section \ref{sez:BIDUAL}. Suppose $r \equiv r_{S,\psi,\eta,\eta'} \in \R^\infty(\Omega)$, and let $f$ be the unique function on $X$ such that ${f \circ \eval[R]} = r$. We want to show that $f \in \C^\infty(X)$. Since $f = {q_{\psi,\eta,\eta'} \circ \eval[S] \circ (\eval[R]|_\Omega)^{-1}}$, it will be enough to show that ${\eval[S] \circ (\eval[R]|_\Omega)^{-1}}$ is a smooth mapping from $X$ into \GL{F}. Note that the groupoid
\begin{equazione}\label{equ:GL(E)xGL(F)}
\GL{E} \xfro{\quad \pr_E \quad} {\GL{E} \times_{M\times M} \GL{F}} \xto{\quad \pr_F \quad} \GL{F}
\end{equazione}
is the product of \GL{E} and \GL{F} in the category of Lie groupoids over $M$; in particular, the projections $\pr_E$ and $\pr_F$ are homomorphisms of Lie groupoids over $M$. We have the following commutative diagram
$$%
\xymatrix@C=50pt{ & \:({R\oplus S})(\Omega)\: \ar[d]_\iso^{\text{homeomorphism}} \ar@{^(->}[r]^-{e_{R,S}} & {\GL{E} \times_{M\times M} \GL{F}} \ar[d]^{\pr_E} \\ \Omega \ar[ur]^(.45){\eval[R\oplus S]} \ar[r]^-{\eval[R]} & \:X = R(\Omega)\: \ar@{^(->}[r]^-{\text{submanifold}} & \GL{E}\text,\!\!}
$$%
where we define $e_{R,S}$ as the map whose composition with the embedding (\ref{equ:GL(E)xGL(F)->GL(E+F)}) equals the inclusion of $({R\oplus S})(\Omega)$ into \GL{E\oplus F}. Now, by Condition (ii), $(\Omega,{R \oplus S})$ is a representative chart. In particular, $({R \oplus S})(\Omega)$ is a submanifold of \GL{E\oplus F} and, consequently, of the product groupoid (\ref{equ:GL(E)xGL(F)}).
Put $\Gamma \equiv {\canhom[\mathcal G]}^{-1}(\Omega)$. By the surjectivity of \canhom[\mathcal G] (Theorem \ref{thm:SURJECT}) and Lemma \ref{lem:Gamma->Sigma}, the maps
$$%
(\varrho \oplus \varsigma)|_\Gamma: \Gamma \longto (R \oplus S)(\Omega) \qquad \text{and} \qquad \varrho|_\Gamma: \Gamma \longto R(\Omega)
$$%
must be open (in fact, submersive) surjections. Then we can apply Corollary \ref{cor:HOMEO=>DIFFEO} with $\psi \equiv \pr_E$, $\mathcal G' \equiv \mathcal G$, $\rho \equiv (\varrho,\varsigma)$ {\large(}that is, the unique Lie groupoid homomorphism of $\mathcal G$ into (\ref{equ:GL(E)xGL(F)}) which composed with (\ref{equ:GL(E)xGL(F)->GL(E+F)}) yields ${\varrho \oplus \varsigma}${\large)}, $\rho' \equiv \varrho$, and $\Gamma' \equiv \Gamma$, in order to conclude that the transition homeomorphism in the preceding diagram is in fact a {\em diffeomorphism.} This implies the smoothness of the map ${\eval[S] \circ (\eval[R]|_\Omega)^{-1}}: X \to \GL{F}$, as it is clear from the commutativity of
$$%
\xymatrix@C=50pt{ & \:(R\oplus S)(\Omega)\: \ar@{^(->}[r]^-{e_{R,S}}_(.4){\text{(smooth)}} & {\GL{E} \times_{M\times M} \GL{F}} \ar[d]^{\pr_F} \\ X \ar[ur]^(.45){\text{(trans.~diffeo.)}^{-1} \quad} \ar[r]^(.57){(\eval[R]|_\Omega)^{-1}} & \:\Omega\: \ar@{^(->}[r]^-{\eval[S]} \ar[u]_{\eval_{R \oplus S}} & \GL{F}\text.\!\!}
$$%

From Condition (i), it follows that $\R^\infty$ makes \bidual{\mathcal G} a smooth manifold and that for each representative chart $(\Omega,R)$ the restriction $\eval[R]|_\Omega$ induces a diffeomorphism of $\Omega$ onto $R(\Omega)$. Finally, it is not difficult to see, by using the remark at the beginning of the proof of Lemma \ref{lem:Gamma->Sigma} once again and the surjectivity of \canhom[\mathcal G], that the source map of the groupoid \bidual{\mathcal G} is a submersion. This completes the proof that \bidual{\mathcal G} is a Lie groupoid.
\end{proof}

\sloppy
Our (provisional) characterization of the property of parareflexivity for proper Lie groupoids is the following improvement on Proposition \ref{prop:CRIT"1}:

\fussy
\begin{proposizione}[Parareflexivity I]\label{prop:CRIT"2}
A proper Lie groupoid is parareflexive if, and only if, the space of arrows of the respective Tannakian bidual can be covered with domains of representative charts.
\end{proposizione}\begin{proof}
We begin by proving that the condition (ii) of Proposition \ref{prop:CRIT"1} is superfluous, i.e., satisfied by any proper Lie groupoid $\mathcal G$.

We need a general remark about submersions first. Suppose we are given a commutative triangle of maps between smooth manifolds
\begin{equazione}\label{equ:rmk-SUBM}
\begin{split}
\xymatrix@C=70pt@R=4pt{ & X \ar@{-->}[dd]^s \\ Y \ar[dr]_{f'} \ar@{->>}[ur]^f & \\ & X'}
\end{split}
\end{equazione}
where $f$ is a submersion onto $X$, $f'$ is smooth, and $s$ is any map making the triangle commute. Then $s$ is smooth. In particular, when also $f'$ is a surjective submersion, $s$ must be a diffeomorphism whenever it is bijective.

Let a representative chart $(\Omega,R)$ be given, in which, let us say, $R = (E,\varrho)$, and let $S = (F,\varsigma)$ be any representation of $\mathcal G$. Put $\Gamma \equiv (\canhom[\mathcal G])^{-1}(\Omega)$. We know that $\varrho$ induces a submersion of $\Gamma$ onto a submanifold $X \equiv R(\Omega)$ of \GL{E}. Then, in (\ref{equ:rmk-SUBM}), take $f$ to be the restriction $\varrho|_\Gamma: \Gamma \to R(\Omega)$, $f'$ to be the map
$$%
(\varrho,\varsigma)|_\Gamma: \Gamma \longto {\GL{E} \times_{M\times M} \GL{F}}\text, \qquad g \mapsto (\varrho(g),\varsigma(g))\text,
$$%
and $s$ to be
$$%
{(\eval[R],\eval[S]) \circ (\eval[R]|_\Omega)^{-1}}: R(\Omega) \longto {\GL{E} \times_{M\times M} \GL{F}}\text.
$$%
By our remark about submersions, $s$ is a {\em smooth} section to the projection $\pr_E$ (\ref{equ:GL(E)xGL(F)}) and hence, in particular, it must be an {\em immersion.} We contend that $s$ is an embedding of manifolds; from this, it will follow that $(R,S)(\Omega) = s(R(\Omega))$ is a submanifold of the product groupoid (\ref{equ:GL(E)xGL(F)}) and hence that $(\Omega,{R\oplus S})$ is a representative chart. Now, for each open subset $\Lambda$ of \GL{E}, we have
$$%
s\bigl(R(\Omega) \cap \Lambda\bigr) = {s(R(\Omega)) \cap \bigl(\Lambda \times_{M\times M} \GL{F}\bigr)}\text.
$$%
Thus, the map $s$ is an {\em open,} one-to-one correspondence between $R(\Omega)$ and the subspace $s(R(\Omega))$ of the product groupoid (\ref{equ:GL(E)xGL(F)}). Our contention is proven.

To finish the proof, we need to show that \canhom[\mathcal G] is a {\em submersion.} (Recall Definition \ref{def:parareflexive}.) By the condition (i), this follows easily from Lemma \ref{lem:Gamma->Sigma} and the remark that for any representative chart $(\Omega,R)$ the bijection $\eval[R]|_\Omega: \Omega \simeq R(\Omega)$ must be a diffeomorphism.
\end{proof}

\sezione{Main criterion for the smoothness of the bidual}\label{sez:MC}
The present section is a continuation of the preceding one. Here, we derive an explicit, ``algebraic'' characterization of parareflexivity for proper Lie groupoids (no longer involving the notion of representative chart), which is to provide a counterpart to Theorem \ref{thm:REFLEX}. We conclude with some examples.

\medskip

Let $\mathcal G$ be a Lie groupoid over a manifold $M$. We say that a submanifold $X \subset M$ is a \textit{slice at~$\,x \in X$}, if the orbit immersion ${\mathcal G x} \into M$ is transversal to $X$ at $x$. A submanifold $S \subset M$ will be called a \textit{slice,} if it is a slice at each point $s \in S$. Note that if $X$ is a submanifold of $M$, and $g \in \mathcal G^X := \s^{-1}(X)$, then $X$ is a slice at $x = \s(g)$ if, and only if, the intersection ${\mathcal G^X \cap \t^{-1}(x')}$, $x' = \t(g)$ is transversal at $g$. From this remark, it follows that for each submanifold $X$, the subset of all points at which $X$ is a slice forms an open subset of $X$. If a submanifold $S$ is a slice, then the intersection ${\s^{-1}(S) \cap \t^{-1}(S)}$ is transversal, so that the restriction $\mathcal G|_S$ is a Lie groupoid over $S$; moreover, ${\mathcal G\cdot S}$ is an invariant open subset of $M$.

For a proof of the following theorem, we refer the reader to the original paper by \mbox{Zung} \cite{Zu06}. We also mention here an earlier paper by \mbox{Weinstein} \cite{We02}, which covers the regular case.

\begin{enunciato*}[N.~T.~Zung]{Local Linearizability Theorem}
Let $\mathcal G$ be a proper Lie groupoid. Let $x$ be a base point which is not moved by the standard action of~$\,\mathcal G$ on its own base.

Then there exists a continuous linear representation $G \to \GL{V}$ of the isotropy group $G \equiv \mathcal G_x$ on a finite dimensional vector space $V$ such that, for some open neighbourhood $U$ of $x$, one can find an isomorphism of Lie groupoids $\mathcal G|_U \iso {G \ltimes V}$ which makes $x$ correspond to zero.\sd
\end{enunciato*}
From this theorem, it follows that any proper Lie groupoid is, locally in the vicinity of each base point, Morita equivalent to an action groupoid associated with some linear representation of a compact Lie group.

\medskip

We proceed to recall the canonical action of an isotropy group on the normal space to the corresponding orbit \cite{We02,Zu06}.

Let $\mathcal G$ be a Lie groupoid over a manifold $M$. Consider the isotropy group $\mathcal G_m$ at a base point $m \in M$. Given any element $g \in \mathcal G_m$, choose an arbitrary local bisection $\sigma: U \to \mathcal G$, ${\t \circ \sigma}: U \isoto U'$ with $\sigma(m)=g$, and then take the tangent map $\T[m]{(\t \circ \sigma)}: \T[m]{M} \to \T[m]{M}$. Since ${\t \circ \sigma}$ maps the orbit $O_m = {\mathcal G\cdot m}$ into itself, this tangent map induces a well-defined endomorphism of the quotient
\begin{equazione}\label{normtang}
\normtang[m]{M}\: \bydef\: \T[m]{M}/\T[m]{O_m}\text.
\end{equazione}
We observe that this endomorphism does not depend on the choice of $\sigma$: to see this, fix any decomposition of the tangent space $\T[g]{\,\mathcal G} \iso {\T[g]{(\mathcal G^m)} \oplus \T[m]{M}}$, so that \T[g]{\,\s} corresponds to the projection onto the second factor, and note that the image under \T[g]{\,\t} of \T[g]{(\mathcal G^m)} is precisely \T[m]{O_m}. We therefore obtain a well-defined (clearly continuous) Lie group action $\mu_m: \mathcal G_m \longto \mathit{GL}\bigl(\normtang[m]{M}\bigr)$. We put
\begin{equazione}\label{sgp-ineffect}
K_m \bydef\: \kernel{\mu_m}
\end{equazione}
and call this (closed, normal) subgroup of $\mathcal G_m$ the \textit{ineffective part} of $\mathcal G_m$.

Our next task will be to prove the following:

\begin{proposizione}\label{pre-MC}
Let $\mathcal G$ be proper. Let $\varrho: \mathcal G \to \GL{E}$ be a representation on a vector bundle $E$ over the base $M$ of~$\,\mathcal G$. For every base point $m \in M$, the two conditions below are equivalent:
\begin{elenco}
\item[(a)]$D_m := \mathrm{Ker}\bigl[\varrho_m: G_m \to \GL{E_m}\bigr]$ is contained in the ineffective part $K_m$ of $G_m$;
\item[(b)]there exists an open neighbourhood $U \subset M$ of $m$ such that the image $X_U := \varrho\bigl(\mathcal G|_U\bigr)$ is a submanifold of \GL{E}.
\end{elenco}
In (b), we can in fact assume $U$ to be invariant i.e.\ ${\mathcal G\cdot U} = U$.
\end{proposizione}

\noindent Notation: We let $X = \varrho(\mathcal G)$ denote the full image, so that $X_U = X \cap \GL{E}|_U$ is an open subset.

By the Local Linearizability Theorem, it is not restrictive to assume that there exists a linear slice $S$ for $\mathcal G$ at the given base point $m$, such that ${\mathcal G\cdot S} = M$. The proposition will then be an immediate consequence of Proposition \ref{MC:linear-case} below and the following couple of lemmas:

\begin{lemma}\label{lem:m.e.-invariance1}
Let $\phi: \mathcal H \to \mathcal G$ be a Morita equivalence of Lie groupoids. Let $\varrho$ be a representation of~$\,\mathcal G$ on a vector bundle $E$. For any base points $n_0$ of~$\,\mathcal H$ and $m_0 = \phi(n_0)$ of~$\,\mathcal G$, the induced isomorphism of isotropy groups $\phi_{n_0}: \mathcal H_{n_0} \xto\iso \mathcal G_{m_0}$ establishes a bijection between the respective ineffective subgroups, and between the kernels of the isotropy representations $(\phi^* \varrho)_{n_0}: \mathcal H_{n_0} \to \mathit{GL}\bigl([\phi^* E]_{n_0}\bigr)$ and $\varrho_{m_0}: \mathcal G_{m_0} \to \GL{E_{m_0}}$.
\end{lemma}

\begin{lemma}\label{lem:m.e.-invariance2}
Under the same assumptions as in the previous lemma, the image $Y = \image{\!(\phi^* \varrho)}$ is a submanifold of \GL{\phi^*E} if and only if $X = \image{\varrho}$ is a submanifold of \GL{E}.
\end{lemma}

\begin{proof}[Proof of Lemma \ref{lem:m.e.-invariance1}]
At the level of base manifolds, the tangent map $\T[n_0]{\phi}: \T[n_0]{N} \to \T[m_0]{M}$ induces a linear isomorphism of \normtang[n_0]{N} onto \normtang[m_0]{M} which $\phi_{n_0}$-intertwines the actions of $\mathcal H_{n_0}$ and $\mathcal G_{m_0}$ on the latter spaces. In particular, the kernels of the two actions are sent into each other by the group isomorphism $\phi_{n_0}$.

On the other hand, the canonical linear isomorphism $[\phi^* E]_{n_0} \simeq E_{\phi(n_0)} = E_{m_0}$ is $\phi_{n_0}$-equivariant with respect to the group actions $(\phi^* \varrho)_{n_0}$ and $\varrho_{m_0}$.
\end{proof}

\begin{proof}[Proof of Lemma \ref{lem:m.e.-invariance2}]
To begin with, observe that it is no loss of generality to assume $\phi$ to be a surjective submersion $\mathcal H \onto \mathcal G$, $N \onto M$. Indeed, as we already pointed out in Section \ref{sez:BIDUAL}, there always exist Morita equivalences $\psi$ and $\chi$ with this property from a third Lie groupoid $\mathcal P$ into $\mathcal G$ and $\mathcal H$ respectively such that $\psi = {\phi \circ \chi}$ up to some natural transformation of Lie groupoid homomorphisms. One then has an isomorphism ${\psi^* \varrho} \iso {\chi^* \phi^* \varrho}$ of representations of $\mathcal P$.

Let $\gamma: \GL{\phi^* E} \to \GL{E}$ be the homomorphism associated with the canonical vector bundle map ${\phi^* E} \to E$ over $\phi: N \to M$. From the surjectivity of $\phi: \mathcal H \to \mathcal G$, it follows that $Y = \gamma^{-1}(X)$. Since $\gamma$ itself must be a surjective submersion, we conclude at once that $X$ is a submanifold if and only if so is $Y$.
\end{proof}

\sottosezione{Study of the linear case}

We let $\mathcal G = {G\ltimes V}$ throughout, where $G \xto\mu \GL{V}$ is a continuous linear action of a compact Lie group $G$ on a finite dimensional real vector space $V$.

Let $K$ denote $\mathrm{Ker}\bigl[G \xto\mu \GL{V}\bigr]$, and, for each $v \in V$, let $G_v \subset G$ denote the isotropy subgroup at $v$. As usual, $K_v \subset G_v$ will denote the ineffective part of $G_v$. Note that $K_0 = K \subset K_v$ for all $v$.

Let $\varrho: {G\ltimes V} \to \GL{\underline{\nC}^k}$ be any Lie groupoid representation, where $\underline{\nC}^k = {V\times \nC^k}$ is the trivial rank $k$ vector bundle over $V$. By canonically identifying $\GL{\underline{\nC}^k}$ to $V\times V\times \GL{k,\nC}$, we write $\varrho$ in components as follows:
\begin{equazione}\label{components-rho}
\varrho(g,v) = \bigl(v,{g\cdot v},{\boldsymbol\varrho}(g,v)\bigr)
\end{equazione}
where ${g\cdot v} = {\mu(g)}(v)$. For each $v \in V$, we let $\varrho_v: G_v \to \GL{k,\nC}$ be the map $g \mapsto {\boldsymbol\varrho}(g,v)$. We put $D_v = \kernel{\varrho_v} \subset G_v$ and $D = D_0 \subset G$.

\begin{lemma}\label{ker-stability}
$D_v = {D \cap G_v}$ for each $v$.
\end{lemma}\begin{proof}
Let $0 \neq v \in V$ be given. Let $L = \langle v\rangle = \{tv | t \in \nR\}$ be the line through $v$. Since the action $G \xto\mu \GL{V}$ on $V$ is \nR-linear, $G_v = G_{tv}$ for all $t \neq 0$. This immediately implies that the restriction
\begin{equazione}\label{isotropy-isotopy}
{G_v\times L} \subset {G\times V} \xto{\:\boldsymbol\varrho\:} \GL{k}
\end{equazione}
is an isotopy through representations $G_v \to \GL{k}$ (note by contrast that $\boldsymbol\varrho$ need not be an isotopy in general). Since all the representations belonging to such an isotopy must have the same kernel, we get $D_v = {D \cap G_v}$.
\end{proof}

\begin{lemma}\label{factor-through-quotient}
Let $N \subset {K \cap D} \subset G$ be a closed subgroup, normal in $G$, and put $\widetilde G = G/N$. Let $\widetilde\mu: \widetilde G \to \GL{V}$ be the induced action. Then there is a (unique) Lie groupoid representation $\widetilde\varrho: {\widetilde G\ltimes V} \to \GL{\underline{\nC}^k}$ such that
\begin{equazione}\label{quotient-repr}
\begin{split}
\xymatrix{{G\ltimes V}\ar[d]^(.4){\pr \times V}\ar[rr]^\varrho & & \GL{\underline{\nC}^k}\text.\!\! \\ {\widetilde G\ltimes V}\ar@{-->}[urr]_{\widetilde\varrho} & &}
\end{split}
\end{equazione}
\end{lemma}\begin{proof}
If a mapping of sets exists making (\ref{quotient-repr}) commute, it will necessarily be smooth because so is $\varrho$ and because ${\pr \times V}$ is a surjective submersion. Moreover, the same map will be a homomorphism of groupoids. Now, it is straightforward to check that the prescription
\begin{equazione}\label{def:quotient-repr}
{\widetilde\varrho(\tilde g,v)}(v,\xi) = \bigl({g\cdot v},{{\boldsymbol\varrho}(g,v) \xi}\bigr) \quad \text{for all} \quad (v,\xi) \in {V\times \nC^k}
\end{equazione}
is well defined---thanks to Lemma \ref{ker-stability}.
\end{proof}

Suppose the image $X = \varrho(\mathcal G) \subset \GL{\underline{\nC}^k}$ is a submanifold. Then we know from Lemma \ref{lem:Gamma->Sigma} that $\mathcal G \xto\varrho X$ is a (surjective) submersion.

Assume for a moment that the group $G$ is connected. Then $X$ is a connected manifold and therefore the dimension of the fibres of the submersion $\varrho$ is constant over $X$. Hence, putting $x_v = \varrho(e,v) \in X$ for each $v \in V$, where $e$ denotes the unit of the group $G$, we see that $\algdim{\!(D \cap G_v)} = \algdim{\varrho^{-1}(x_v)} = \algdim{\varrho^{-1}(x_0)} = \algdim{D}$ for all $v$. Since $D \cap G_v$ is a closed subgroup of $D$, this implies that the identity component $D^{(e)}$ of $D$ is contained in $G_v$. Thus
\begin{equazione}\label{NEC1}
D^{(e)} \subset {\bigcap_{v \in V} G_v} = K\text.
\end{equazione}
For a nonconnected $G$, it will be enough to consider the restriction $\varrho^{(e)}$ of $\varrho$ to ${G^{(e)} \ltimes V}$, where $G^{(e)}$ denotes the identity component of $G$, and observe that the condition (\ref{NEC1}) for $\varrho$ is implied by the analogous condition for $\varrho^{(e)}$.

The foregoing reasonings show that (\ref{NEC1}) is a necessary condition in order for the image of $\varrho$ to be a submanifold of \GL{\underline{\nC}^k}. However, this condition is not sufficient: take for example $G = \nZ/2 = \{\pm1\}$ acting on $V = \nR$ by scalar multiplication, and $\varrho$ trivial.

Consider the intersection ${K \cap D} \subset G$. This is a closed normal subgroup of $G$ and hence we can apply Lemma \ref{factor-through-quotient} with $N = {K \cap D}$. Note that, in view of (\ref{NEC1}), the image $\widetilde D$ of $D$ under the quotient map $G \to \widetilde G = G/(K \cap D)$ is a finite subgroup of $\widetilde G$. Moreover, one has ${\widetilde K \cap \widetilde D} = \{1\}$ where $\widetilde K$ is the image of $K$ in $\widetilde G$. Since $\image{\varrho} = \image{\widetilde\varrho}$ in (\ref{quotient-repr}), we may accordingly assume\inciso{without loss of generality}the following to hold:
\begin{equazione}\label{ASS1}
D\text{~finite,} \qquad D\cap K = \{e\}\text.
\end{equazione}

The finiteness assumption on the kernel $D$ of the isotropy representation $\varrho_0: G \to \GL{k}$ implies at once the injectivity of the Lie algebra map $\T[e]{\varrho_0}: \T[e]{G} \to \End{\nC^k}$. The matrix of the differential of the mapping (\ref{components-rho}) of ${G\times V}$ into $V\times V\times \GL{k}$ at the point $(e,0)$ is of the form
\begin{equazione}\label{matrix-diff-comp-rho}
\begin{pmatrix}
\id_V & 0 \\ * & * \\ * & \T[e]{\varrho_0}
\end{pmatrix}\text.
\end{equazione}
Since $\varrho: \mathcal G \to X$ is a submersion, the injectivity of the differential ${\D\varrho\,}(e,0): \T[(e,0)]{\mathcal G} \to \T[x_0]{X}$ shows that $\varrho$ is a local diffeomorphism at $(e,0)$. We contend that the same is true at each point $(d,0)$ with $d \in D$. Now, $(g,v) \mapsto ({d\cdot g},v)$ defines a diffeomorphism of ${G\times V}$ into itself, mapping $(e,0)$ to $(d,0)$. Similarly,
\begin{equazione}\label{self-diffeo-auxilary}
(v,v',\lambda) \mapsto \bigl(v,{d\cdot v'},{\boldsymbol\varrho(d,v') \circ \lambda}\bigr)
\end{equazione}
defines an automorphism of $V\times V\times \GL{k}$ which leaves the submanifold $X$ invariant and fixes the point $x_0 = \varrho(e,0) = (0,0,\mathbf I_k) \in X$. The mapping (\ref{components-rho}) intertwines these two automorphisms. Our claim follows at once by taking the corresponding tangent maps.

For each $d \in D$, we choose an open neighbourhood $\Gamma_d$ of $(d,0)$ in $\mathcal G$ such that $\varrho$ restricts to a diffeomorphism between $\Gamma_d$ and an open neighbourhood $U_d = \varrho(\Gamma_d)$ of $x_0$ in $X$. We may assume the $\Gamma_d$ to be pairwise disjoint. If we put $U = \textcap{d\in D}{}{U_d}$, we see that for each $x \in U$ the fibre $\varrho^{-1}(x) \subset \mathcal G$ contains at least \norm{D} distinct points. On the other hand, for each $v \in V$ the fibre $\varrho^{-1}(x_v) = {D \cap G_v}$ can at most contain \norm{D} points. Since the map $v \mapsto x_v$ of $V$ into $X$ is continuous, and since $U$ is a neighbourhood of $x_0$, it follows that $D \subset G_v$ for all $v$ sufficiently close to zero and therefore for all $v$ (recall that $G_{tv} = G_v$ for $t \neq 0$). Thus,
$$%
D \subset {D \cap {\bigcap_{v\in V} G_v}} = {D \cap K} = \{e\}\text.
$$%

Summing up, we have shown\inciso{now in full generality}that a necessary condition for the image of $\varrho$ to be a submanifold of \GL{\underline\nC^k} is that
\begin{equazione}\label{NEC2}
D \subset K\text.
\end{equazione}
This condition proves to be also sufficient. Indeed, if the inclusion (\ref{NEC2}) holds, we can apply Lemma \ref{factor-through-quotient} with $N = D$ and then observe that the representation $\widetilde\varrho$ is globally faithful and hence maps ${\widetilde G\ltimes V}$ onto a submanifold of \GL{\underline\nC^k}, in view of the following lemma:

\begin{lemma}\label{lem:FAITH=>SUBMAN}
Let $\mathcal G$ be a proper Lie groupoid, and let $\varrho: \mathcal G \to \GL{E}$ be a representation. Suppose that $\varrho$ is faithful. Then the image $\varrho(\mathcal G)$ is a submanifold of \GL{E}.
\end{lemma}\begin{proof}
To begin with, we observe that for any given arrow $g$ in $\mathcal G$, and for each neighbourhood $\Gamma$ of $g$, there is an open neighbourhood $P$ of $\varrho(g)$ in \GL{E} such that $\varrho^{-1}(P) \subset \Gamma$. This can be seen as in the proof of Theorem \ref{thm:INJ=>HOMEO}.

By Lemma \ref{lem:INJREP}, each arrow $g$ admits an open neighbourhood $\Gamma_g$ such that $\varrho$ induces a smooth isomorphism between $\Gamma_g$ and a submanifold of $\GL{E}$. One can then choose an open neighbourhood $P_g \subset \GL{E}$ of $\varrho(g)$ such that $\varrho^{-1}(P_g) \subset \Gamma_g$. For any given pair $x,x'$ of base points, put $\Gamma_{x,x'} \equiv {\bigcup \varrho^{-1}(P_g)}$, the union being taken over all $g \in \mathcal G(x,x')$. We claim that $\varrho$ induces a smooth isomorphism between $\Gamma_{x,x'}$ and a submanifold of $\GL{E}$. By construction, $\varrho$ restricts to an immersion of $\Gamma_{x,x'}$ into $\GL{E}$. For each $g \in \mathcal G(x,x')$,
$$%
{\varrho(\Gamma_{x,x'}) \cap P_g} = \varrho\bigl(\varrho^{-1}(P_g)\bigr)
$$%
is an open subset of the submanifold $\varrho(\Gamma_g) \subset \GL{E}$, because $\varrho$ is a smooth isomorphism of $\Gamma_g$ onto $\varrho(\Gamma_g)$. Since the open sets $P_g$ cover $\varrho(\Gamma_{x,x'})$ as $g$ ranges over $\mathcal G(x,x')$, $\varrho(\Gamma_{x,x'})$ is a submanifold of $\GL{E}$. Moreover, since $\varrho$ is a local smooth isomorphism of $\Gamma_{x,x'}$ onto $\varrho(\Gamma_{x,x'})$, it will be a global diffeomorphism provided it is injective on $\Gamma_{x,x'}$. Now, $\varrho(\gamma') = \varrho(\gamma)$ implies $\gamma', \gamma \in \varrho^{-1}(P_g) \subset \Gamma_g$ for at least one $g$ and therefore $\gamma' = \gamma$, because $\varrho$ is injective over $\Gamma_g$.

Finally, one application of the usual properness argument will yield open neighbourhoods $B \ni x$ and $B' \ni x'$ in $M$ such that the set $\mathcal G(B,B')$ is contained in $\Gamma_{x,x'}$. This essentially finishes the proof.
\end{proof}

Putting everything together, we conclude

\begin{proposizione}\label{MC:linear-case}
A representation $\varrho: {G\ltimes V} \to \GL{\underline\nC^k}$ maps ${G\ltimes V}$ onto a submanifold of \GL{\underline\nC^k} if, and only if, the kernel of its isotropy representation at zero acts ineffectively on $V$.
\end{proposizione}

\begin{paragrafetto}[Remark]\label{main-condition=>local}
The condition (\ref{NEC2}) implies $D_v \subset K_v$ for all $v \in V$ (notation introduced before Lemma \ref{ker-stability}). Indeed, by Lemma \ref{ker-stability}, if $D \subset K$ then
$$%
D_v = {D \cap G_v} \subset {K \cap G_v} \subset K_v\text.
$$%
\end{paragrafetto}

\sottosezione{The Main Criterion and its application to the regular case}

A representative chart of the form $\bigl(\bidual{\mathcal G}|_U,R\bigr)$ will be said to be \textit{diagonal.}

\begin{lemma}\label{diag-repr-charts}
For $\mathcal G$ proper, \bidual{\mathcal G} is smooth if, and only if, for each base point $m$ there exists a diagonal representative chart at $m$.
\end{lemma}\begin{proof}
(Sufficiency.)~~Let $\lambda \in \bidual{\mathcal G}(m,m')$ be given. We want to construct a representative chart $(\Omega,R)$ in a neighbourhood $\Omega$ of $\lambda$. Choose $g \in \mathcal G(m,m')$ with $\canhom(g) = \lambda$ (Theorem \ref{thm:SURJECT}) and a local bisection $\sigma: U \into \mathcal G$, ${\t \circ \sigma}: U \isoto U'$ with $\sigma(m)=g$. It is not restrictive to assume there exists a diagonal representative chart $\bigl(\bidual{\mathcal G}|_U,R\bigr)$. Then we can take $\Omega = \bidual{\mathcal G}(U,U')$, for we have the following commutative diagram:
$$%
\xymatrix{{\mathcal G(U,U')}\ar[r]^-\canhom & \bidual{\mathcal G}(U,U')\ar[r]^-{\eval[R]} & {\GL{E}}(U,U') \\ {\mathcal G|_U}\ar[r]^-\canhom\ar[u]^\iso_{\sigma\text-} & \bidual{\mathcal G}|_U\ar[r]^-{\eval[R]}\ar[u]^\iso_{(\canhom\circ\sigma)\text-} & {\GL{E}}|_U\ar[u]^\iso_{(\varrho\circ\sigma)\text-}}
$$%
(where $R = (E,\varrho)$, and where $\sigma\text-$ etc.\ stand for left multiplication in the appropriate sense).

(Necessity.)~~By a compactness argument, we can find a finite set of representative charts $(\Omega_i,R_i)$ covering the isotropy group $\bidual{\mathcal G}|_m$. Put $R \equiv {\smash{\bigoplus}_i R_i}$, where, let us say, $R = (E,\varrho)$. Since \bidual{\mathcal G} is smooth by assumption, it follows that each $(\Omega_i,R)$ is a representative chart as well, by Proposition \ref{prop:CRIT"1}. By taking, if necessary, a larger $R$, we can also assume that the evaluation representation \eval[R] is faithful on $\bidual{\mathcal G}|_m$. As in the proof of Lemma \ref{lem:FAITH=>SUBMAN}, one can then show that for each $\lambda \in {\Omega_i \cap \bidual{\mathcal G}|_m}$ there exists some open ball $P_\lambda \subset \GL{E}$ centred at $\eval[R](\lambda)$ such that ${\eval[R]}^{-1}(P_\lambda) \subset \Omega_i$. Put $\Omega \equiv {\smash{\bigcup}_\lambda {\eval[R]}^{-1}(P_\lambda)}$. Again as in the above-mentioned proof, one can check that $(\Omega,R)$ is a representative chart. This evidently finishes the proof of the lemma.
\end{proof}

\begin{teorema}[Main Criterion for smoothness]\label{MC}
Let $\mathcal G$ be a proper Lie groupoid. In order that also the Tannakian bidual \bidual{\mathcal G} may be a Lie groupoid, it is necessary and sufficient that for each base point $m$ there exists a representation $\varrho: \mathcal G \to \GL{E}$ whose associated isotropy representation at $m$ has ineffective kernel:
\begin{equazione}\label{MC-cond}
\mathrm{Ker}\bigl[\varrho_m: \mathcal G|_m \to \GL{E_m}\bigr] \subset K_m\text.
\end{equazione}
\end{teorema}\begin{proof}
(Sufficiency.)~~Let us fix a representation $R = (E,\varrho)$ such that \eval[R] injects $\bidual{\mathcal G}|_m$ into \GL{E}. Of course, then $\kernel{\varrho_m} \subset K_m$ (by our assumption) and therefore, by Proposition \ref{pre-MC}, there exists an open neighbourhood $U$ of $m$ such that $\eval[R]\bigl(\bidual{\mathcal G}|_U\bigr) = \varrho\bigl(\mathcal G|_U\bigr)$ is a submanifold of \GL{E}.

All we still lack in order to get a diagonal representative chart at $m$ is the injectivity of the map \eval[R] on a (possibly smaller) neighbourhood of $\bidual{\mathcal G}|_m$ of the same form. It is clearly enough to show that we can find a linear slice $S$ at $m$ such that \eval[R] is injective on each isotropy group in $\bidual{\mathcal G}|_S$. We may therefore assume that $\varrho$ is, for some compact Lie group $G$, a representation of some linear action groupoid ${G\ltimes V}$ on a (trivial) vector bundle over $V$. Any other representation of the groupoid $\mathcal G$ will induce a representation $\varrho'$ of ${G\ltimes V}$ with $\kernel{{\varrho'}_0} \supset \kernel{\varrho_0}$. If $\kernel{\varrho_v} \nsubset \kernel{{\varrho'}_v}$ were true for some $v \in V$ then, by Lemma \ref{ker-stability}, $\kernel{\varrho_0} \nsubset \kernel{{\varrho'}_0}$ would also be true: contradiction. Hence $\kernel{{\varrho'}_v} \supset \kernel{\varrho_v}$ for all $v$. This shows the required injectivity.

(Necessity.)~~Let $\bigl(\bidual{\mathcal G}|_U,R\bigr)$ be an arbitrary diagonal representative chart at $m$, and let $R = (E,\varrho)$. Then, by definition, $\varrho\bigl(\mathcal G|_U\bigr) = \eval[R]\bigl(\bidual{\mathcal G}|_U\bigr)$ is a submanifold of \GL{E} and, therefore, (\ref{MC-cond}) follows from Proposition \ref{pre-MC}.
\end{proof}

Let us call a representation $\varrho: \mathcal G \to \GL{E}$ \textit{effective at $m$} when the inclusion (\ref{MC-cond}) holds. We shall say that \textit{$\mathcal G$ has enough effective representations} if the condition in the statement of Theorem \ref{MC} is satisfied. Then our final characterization of parareflexivity for proper Lie groupoids reads:

\begin{teorema}[Parareflexivity II]\label{thm:PARAREFLEX}
A proper Lie groupoid is parareflexive if and only if it possesses enough effective representations.\qed
\end{teorema}

We shall now give an immediate application of the preceding criterion. Recall that a Lie groupoid $\mathcal G$ over a manifold $M$ is said to be \textit{regular} when the anchor map $\rho: \mathfrak{g} \to \T{M}$ of the Lie algebroid of $\mathcal G$ \cite{MoeMrc03} has locally constant rank as a morphism of vector bundles over $M$. If $\mathcal G$ is regular then the image of the anchor map $\rho$ is a subbundle of the tangent bundle of $M$ which is also integrable and hence determines a foliation of $M$. By a result of \mbox{Moerdijk} \cite{Moe03}, any regular Lie groupoid $\mathcal G$ over $M$ fits into a short exact sequence of Lie groupoid homomorphisms over $M$
\begin{equazione}\label{sh.ex.sequ1}
1 \to \mathcal B \into \mathcal G \onto \mathcal F \to 1
\end{equazione}
where $\mathcal B$ is a bundle of connected Lie groups, $\mathcal F$ is a foliation groupoid and $\onto$ is a submersion with connected fibres. Observe that, being $\onto$ identical on the base $M$, the epimorphism of isotropy groups $\mathcal G_m \onto \mathcal F_m$ that $\onto$ induces at any base point $m \in M$ will carry the canonical action of $\mathcal F_m$ on the normal space (\ref{normtang}) to the canonical action of $\mathcal G_m$ on the same space; in particular, an arrow in $\mathcal G_m$ will be ineffective if, and only if, so is its image in $\mathcal F_m$.

\begin{corollario}\label{regular+proper=>smooth}
For an arbitrary regular proper Lie groupoid $\mathcal G$, the Tannakian bidual \bidual{\mathcal G} is a Lie groupoid.
\end{corollario}\begin{proof}
We claim there exists a representation of $\mathcal G$ whose kernel is precisely the ineffective part of the isotropy of $\mathcal G$. By the foregoing remarks, we can assume that $\mathcal G$ is a foliation groupoid.

Recall that any foliation groupoid is \mbox{Morita} equivalent to an \'etale groupoid \cite{MoeMrc03,CraMoe01}. On the other hand, the condition (\ref{MC-cond}) in our criterion\inciso{more correctly, the property of existence of a representation satisfying that condition at each base point}is immediately seen to be stable under passage to a Morita equivalent groupoid (Lemma \ref{lem:m.e.-invariance1}). Also recall that if $\mathcal E$ is an arbitrary \'etale groupoid then there is a canonical short exact sequence of groupoids over the base of $\mathcal E$
\begin{equazione}\label{sh.ex.sequ2}
1 \to \mathcal K \into \mathcal E \onto \widetilde{\mathcal E} \to 1
\end{equazione}
in which $\widetilde{\mathcal E}$ is an effective \'etale groupoid and $\mathcal K$ is precisely the ineffective part of $\mathcal E$. The groupoid $\widetilde{\mathcal E}$ is called the \textit{effect} of $\mathcal E$ (\cite{MoeMrc03} p.\ 136). Now, it only remains to observe that each proper effective \'etale groupoid possesses a canonical faithful representation on the tangent bundle to the base manifold, namely the \textit{tangent representation.} For details, see \cite{2008}, Prop.\ 28.4.
\end{proof}

\begin{paragrafetto}[Examples in low orbit codimensions]\label{LowOrbCodim}
We mention, without proof, the following result (to appear). Let $\mathcal G$ be any Lie groupoid over a manifold $M$. We define the \textit{orbit codimension of~$\,\mathcal G$} to be the integer
\begin{equazione}\label{OrbCodim}
\algdim{M/\mathcal G} := {\sup_{x \in M} \: \mathrm{dim}_\nR\bigl(\normtang[x]{M}\bigr)}\text,
\end{equazione}
the notation being as in (\ref{normtang}). One can then show that for any proper Lie groupoid $\mathcal G$, {\em $\algdim{M/\mathcal G} \leqq 2$ implies $\mathcal G$ parareflexive} (the only nontrivial case here is when $\algdim{M/\mathcal G} = 2$). In particular, any proper Lie groupoid over a two dimensional base manifold must be parareflexive.\sd
\end{paragrafetto}

\begin{paragrafetto}[Open problem]\label{OpenProb}
So far, we have not been able to find any example of a proper Lie groupoid which is not parareflexive. We know (from \ref{LowOrbCodim}) that for any such groupoid, the dimension of the corresponding base manifold must be at least three.
\end{paragrafetto}

\appendix

\hbadness=10000
\bibliography{../../shared/biblio/lgpd,../../shared/biblio/tnkd}

\end{document}